\documentclass[11pt]{amsart}
\usepackage{amssymb}
\usepackage{amsthm}

\newtheorem{theorem}{Theorem}[section]
\newtheorem*{thm}{Theorem}
\newtheorem{corollary}[theorem]{Corollary}
\newtheorem{lemma}[theorem]{Lemma}
\newtheorem{proposition}[theorem]{Proposition}
\theoremstyle{definition}
\newtheorem{definition}[theorem]{Definition}

\newtheorem{question}[theorem]{Question}
\theoremstyle{remark}
\newtheorem{remark}[theorem]{Remark}
\newtheorem{example}[theorem]{Example}

\newtheorem{problem}{\bf Problem}
\numberwithin{equation}{section}

\def\R {{\mathbb{R}}}
\def\T {{\mathbb{T}}}

\def\N{{\mathbb{N}}}

\def\i{{\rm i}}

\def\Z {{\mathbb{Z}}}

\def\3{{|\!|\!|}}

\begin{document}

\title{On equivalence relations induced by Polish groups}
\author{Longyun Ding}
\address{School of Mathematical Sciences and LPMC, Nankai University, Tianjin, 300071, P.R.China}
\email{dingly@nankai.edu.cn}
\thanks{Research is partially supported by the National Natural Science Foundation of China (Grant No. 11725103).}
\author{Yang Zheng}
\address{Academy of Mathematics and Systems Science, Chinese Academy of Sciences, East Zhong Guan Cun Road No. 55, Beijing 100190, P.R.China}
\email{yangz@amss.ac.cn}



\begin{abstract}
The motivation of this article is to introduce a kind of orbit equivalence relations which can well describe structures and properties of Polish groups from the perspective of Borel reducibility. Given a Polish group $G$, let $E(G)$ be the right coset equivalence relation $G^\omega/c(G)$, where $c(G)$ is the group of all convergent sequences in $G$.

Let $G$ be a Polish group. (1) $G$ is a discrete countable group containing at least two elements iff $E(G)\sim_BE_0$;
(2) if $G$ is TSI uncountable non-archimedean, then $E(G)\sim_BE_0^\omega$; (3) $G$ is non-archimedean iff $E(G)\le_B=^+$;
(4) if $H$ is a CLI Polish group but $G$ is not, then $E(G)\not\le_BE(H)$; (5) if $H$ is a non-archimedean Polish group but $G$ is not, then $E(G)\not\le_BE(H)$.

The notion of $\alpha$-l.m.-unbalanced Polish group for $\alpha<\omega_1$ is introduced. Let $G,H$ be Polish groups, $0<\alpha<\omega_1$. If $G$ is $\alpha$-l.m.-unbalanced but $H$ is not, then $E(G)\not\le_B E(H)$.

For TSI Polish groups, the existence of Borel reduction is transformed into the existence of a well-behaved continuous mapping between topological groups.
As its applications, for any Polish group $G$, let $G_0$ be the connected component of the identity element $1_G$. Let $G$ and $H$ be two separable TSI Lie groups.
If $E(G)\le_BE(H)$, then there exists a continuous locally injective map $S:G_0\to H_0$. Moreover, if $G_0,H_0$ are abelian, $S$ is a group homomorphism.
In particular, for $c_0,e_0,c_1,e_1\in\N$, $E(\R^{c_0}\times\T^{e_0})\le_BE(\R^{c_1}\times\T^{e_1})$ iff $e_0\le e_1$ and $c_0+e_0\le c_1+e_1$.
\end{abstract}
\maketitle

\section{Introduction}

In recent years, a theory has been developed abundantly in descriptive set theory: using Borel reducibility to investigate the complexity of equivalence relations. Equivalence relations from different branches of mathematics have been systematically studied, and a large number of results on Borel reducibility or non-reducibility have been found. Polish groups are important tools in the study of equivalence relations. The {orbit equivalence relations} generated by the continuous~(or Borel) group actions of Polish groups $G$ on Polish spaces $X$~(denoted by $E_G^X$) account for most of the various equivalence relations concerned by various branches of classical mathematics.

We let $E$ be a Borel equivalence relation on a Polish space $X$. We say that $E$ is {\it smooth} if $E$ is Borel reducible to the equality relation of some Polish space; and say that $E$ is {\it countable} if every $E$-class is countable. We call $E$ {\it essentially countable} if there is a countable
Borel equivalence relation $F$ such that $E$ is Borel reducible to $F$. Feldman-Moore Theorem~\cite[Theorem 7.1.4]{gaobook} asserts that, for every countable Borel equivalence relation $E$ on a Polish space $X$, there is a countable group $\Gamma$ and a Borel action of $\Gamma$ on $X$ such that $E=E^X_\Gamma$.

A very prevalent phenomenon is that the structural properties of a Polish group $G$ can affect the Borel complexity of the orbit equivalence relations $E_G^X$.
For instances, Solecki ~\cite{Solecki} proved that a Polish group $G$ is {compact} if and only if $E_G^X$ is  smooth for any Polish $G$-space $X$;
Kechris~\cite{Kechris_locallycom} proved that if $G$ is a {locally compact Polish group}, then $E_G^X$ is essentially countable for any Polish $G$-space $X$;
Gao and Jackson~\cite{GaoJackson} proved that if $G$ is a {countable discrete abelian group}, then $E_G^X$ is Borel reducible to $E_0$ for any Polish $G$-space $X$;
Hjorth's turbulence theory~\cite{hjorth} asserts that orbit equivalence relations induced by turbulent actions are not Borel reducible to $E_{S_\infty}^X$ for any Polish $S_\infty$-space $X$;
the orbit equivalence relation generated by a continuous action of a {CLI} Polish group on a Polish space $X$ is {pinned}, which implies that any non-pinned Borel equivalence relation is not Borel reducible to $E_G^X$~(see~\cite[\S 17]{kanovei}).
Naturally, we ask the following problem:

\begin{problem}\label{probbb}
Can we characterize  structures and properties of Polish groups through the study of equivalence relations?
\end{problem}

Theoretically, this should be possible: for each kind of Polish groups, studying the properties of all/part of their orbit equivalence relations will inevitably reflect the properties of this kind of Polish groups. However, systematic research needs to find a specific orbit equivalence relation that can well reflect the properties of a given Polish group.

Given a countable discrete group $\Gamma$. There is a
canonical continuous {\it shift action} of $\Gamma$ on $\{0,1\}^\Gamma$ given by $(gp)(h)=p(g^{-1}h)$ for $p\in \{0,1\}^\Gamma$ and $g,h\in \Gamma$.
The {\it free part} of the action is the set
$$(2)^\Gamma=\{p\in\{0,1\}^\Gamma:gp\neq p\mbox{ for all } 1_\Gamma\neq g\in \Gamma\}.$$
We let $E_\Gamma$ be the corresponding orbit equivalence relation on $(2)^\Gamma$. Indeed, the Borel complexity of the relation $E_\Gamma$ can reflect the properties of the group $\Gamma$. For instance, if $E_\Gamma$ is Borel reducible to $E_0$, then $\Gamma$ is amenable, i.e., there is a left-invariant finitely additive probability measure $\mu$ on $\Gamma$~(cf.~\cite[\S 7.4]{gaobook}). For a special case of $\Gamma$, Thomas proved that $E_\Gamma$ is Borel reducible to $E_\Lambda$ if and only if $\Gamma$ embeds into $\Lambda$~(cf.~\cite[Corollary 3.8]{Thomas}). So in this special case, we can think that the equivalence relation $E_\Gamma$ captures the algebraic structure of $\Gamma$ very well.

The main purpose of this article is to attempt to provide a possible answer to Problem~\ref{probbb} for general Polish groups.

Given a Polish group $G$, we define an equivalence relation $E(G)$ on $G^\omega$ as: for $x,y\in G^\omega$,
$$xE(G)y\iff\lim_nx(n)y(n)^{-1}\mbox{ converges in }G.$$
We say $E(G)$ is the \textit{equivalence relation induced by $G$}. Indeed, $E(G)$ is the orbit equivalence relation induced by the action of left multiplication of $c(G)$ on $G^\omega$, where $c(G)$ is the Polish group consisting of all convergent sequences in $G$.

It is clear that all equivalence relations induced by Polish groups are ${\mathbf\Pi}^0_3$.

The study of equivalence relations induced by Polish groups can distinguish Polish groups very well. For instance, we can compare them with some benchmark equivalence relations: $E_0,E_0^\omega$ and $=^+$ (the definitions of these benchmark equivalence relations are left for the next section).

\begin{theorem}
Let $G$ be a Polish group, then we have
\begin{enumerate}
\item[(1)] $G$ is a discrete countable group containing at least two elements iff $E(G)\sim_BE_0$;
\item[(2)] if $G$ is TSI uncountable non-archimedean, then $E(G)\sim_BE_0^\omega$;
\item[(3)] $G$ is non-archimedean iff $E(G)\le_B=^+$.
\end{enumerate}
\end{theorem}

Recently, the authors~\cite{DZtsi} use the Borel complexity of $E(G)$ to characterize TSI non-archimedean Polish groups.

\begin{thm}[{\cite[Theorem 1.3]{DZtsi}}]
Let $G$ be a Polish group. Then the following are equivalent:
\begin{enumerate}
  \item $G$ is TSI non-archimedean;
  \item $E(G)\leq_B E_0^\omega$; and
  \item $E(G)\leq_B \R^\omega/c_0$.
\end{enumerate}
In particular, $E(G)\sim_B E_0^\omega$ iff $G$ is TSI uncountable non-archimedean.
\end{thm}

Many notions about Borel reductions can be applied in the research on equivalence relations induced by Polish groups. For instance, using the notion of right $\iota$-embedability introduced by Lupini and Panagiotopoulos~\cite{LP}, we get the following:

\begin{theorem}
Let $G$ be a non-CLI Polish group and $H$ a CLI Polish group. Then $E(G)\not\le_BE_H^Y$ for any Polish $H$-space $Y$.

In particular, since $c(H)$ is also CLI, we have $E(G)\not\le_B E(H)$.
\end{theorem}

Similarly, using the notion of turbulence introduced by Hjorth~\cite{hjorth} and $c_0$-equality introduced by Farah~\cite{farah}, we get the following:

\begin{theorem}
Let $G$ and $H$ be two Polish groups. If $H$ is non-archimedean but $G$ is not, then $E(G)\not\le_BE_H^Y$ for any Polish $H$-space $Y$.

In particular, since $c(H)$ is also non-archimedean, we have $E(G)\not\le_B E(H)$.
\end{theorem}

Recently, similar result concerning TSI Polish groups is obtained by the authors~\cite{DZtsi}.

\begin{thm}[{\cite[Theorem 1.2]{DZtsi}}]
Let $G,H$ be two Polish groups. If $H$ is TSI but $G$ is not, then $E(G)\not\le_BE(H)$.
\end{thm}

More attractive results appear in the study of CLI Polish groups. Using the notion of $\alpha$-unbalanced relations introduced by Allison and Panagiotopoulos~\cite{AP}, we define the notion of $\alpha$-l.m.-unbalanced Polish groups for $\alpha<\omega_1$ and obtain the following:

\begin{theorem}
Let $G,H$ be Polish groups, $0<\alpha<\omega_1$. If $G$ is $\alpha$-l.m.-unbalanced but $H$ is not, then $E(G)\not\le_B E(H)$.
\end{theorem}

We will also present examples of $\alpha$-l.m.-unbalanced groups. Recall that a topological group $G$ is said to be distal if $1_G\notin\overline{\{ghg^{-1}:g\in G\}}$ for all $h\ne 1_G\in G$. It is known that all TSI Polish groups and all nilpotent Polish groups are distal (c.f.~\cite{rosenblatt}).

\begin{theorem}
Let $G$ be a Polish group.
\begin{enumerate}
\item[(1)] If $G$ is not distal, then $G$ is $1$-l.m.-unbalanced. In particular, if $G$ is locally compact, then $G$ is not distal iff it is $1$-l.m.-unbalanced.
\item[(2)] If $G$ is locally compact, then $G$ is not $2$-l.m.-unbalanced.
\item[(3)] Let $\Lambda$ be an infinite countable discrete group, and $\alpha<\omega_1$. Then $G$ is $\alpha$-l.m.-unbalanced iff the wreath product $\Lambda\wr G$ is $(\alpha+1)$-l.m.-unbalanced.
\end{enumerate}
\end{theorem}

Therefore, we can find a sequence $(G_n)$ of CLI Polish groups such that
$$E(G_0)<_B\cdots<_BE(G_n)<_BE(G_{n+1})<_B\cdots.$$
Unfortunately, we have not found any $\omega$-l.m.-unbalanced CLI Polish group so far, though we know that $S_\infty$ is $\alpha$-l.m.-unbalanced for any ordinal $\alpha$.

The most refined result appears in the study of TSI Polish groups. Applying the tool of additive reduction and ultrafilter limit, the existence of Borel reduction is transformed into the existence of a well-behaved continuous mapping between Polish groups.

Let $G$ be a Polish group. We define equivalence relation $E_*(G)$ on $G^\omega$ as:
$$xE_*(G)y\iff\lim_nx(0)x(1)\cdots x(n)y(n)^{-1}\cdots y(1)^{-1}y(0)^{-1}\mbox{ converges}$$
for $x,y\in G^\omega$. It is clear that $E(G)\sim_BE_*(G)$~(cf.~\cite[Proposition 2.2]{DZ}). It turns out that, for TSI Polish groups, $E_*(G)$ is a more convenient research object than $E(G)$.

\begin{theorem}[Pre-rigid Theorem]
Let $G,H$ be two TSI Polish groups. Suppose $G$ is connected, $H$ is locally compact and the interval $[0,1]$ embeds in $H$. Then $E(G)\le_BE(H)$ iff there exists a continuous map $S:G\to H$ with $S(1_G)=1_H$ such that, for $x,y\in G^\omega$, if $\lim_nd_G(x(n),y(n))=0$, then
$$xE_*(G)y\iff S^\omega(x)E_*(H)S^\omega(y),$$
where the map $S^\omega:X^\omega\to Y^\omega$ is defined as: for $x\in X^\omega$ and $n\in\omega$,
$$S^\omega(x)(n)=S(x(n)).$$
\end{theorem}

To illustrate the strength of this theorem, we recall some notions. For any Polish group $G$, let $G_0$ be the connected component of the identity element $1_G$, which is a closed normal subgroup of $G$. It is known that, if $G$ is locally compact TSI, then $G_0$ is the product of a compact connected Polish group and a group $\R^n$. Now we apply Pre-rigid Theorem on Lie groups as follows.

\begin{theorem}
Let $G$ and $H$ be two separable TSI Lie groups. If $E(G)\le_BE(H)$, then there exists a continuous locally injective map $S:G_0\to H_0$.

In particular, if $G_0,H_0$ are abelian, then $S$ is a group homomorphism.
\end{theorem}

As a corollary, we have $\dim(G)\le\dim(H)$. Furthermore, if $\dim(G)=\dim(H)$ and $G_0$ is compact, then $H_0$ is also compact and $G_0$ is a covering space of $H_0$, and hence the corresponding Lie algebras of $G_0$ and $H_0$ are isomorphic. In particular, we obtain a complete answer for abelian connected Lie groups. Recall that $\T$ is the multiplicative group of all complex numbers with modulus $1$.

\begin{theorem}
Let $c_0,e_0,c_1,e_1\in\N$, then $E(\R^{c_0}\times\T^{e_0})\le_BE(\R^{c_1}\times\T^{e_1})$ iff $e_0\le e_1$ and $c_0+e_0\le c_1+e_1$.
\end{theorem}

Recently, the aforementioned results have been vigorously generalized by the authors~\cite{DZtsi}.
A completely metrizable topological group $G$ is called a \textit{pro-Lie group} if every open neighborhood of $1_G$ contains a normal subgroup $N$ such that $G/N$ is a Lie group (cf.~\cite[Definition 1]{HM07}). By~\cite[Theorem 3.6]{HMS}, every locally compact TSI Polish group is a pro-Lie group.
The following is a Rigid Theorem concerning locally compact connected TSI Polish groups and TSI pro-Lie groups.

\begin{thm}[{\cite[Theorem 6.9]{DZtsi}}]
Let $G$ be a locally compact connected TSI Polish group, $H$ a TSI pro-Lie Polish group. Then $E(G)\leq_B E(H)$ iff there exists a continuous homomorphism $S:G\rightarrow H$ such that $\ker(S)$ is non-archimedean.
\end{thm}

We also get a few results for disconnected Lie groups. Let ${\rm Inn}_G(G_0)$ be the group of all $\iota_u:G_0\to G_0$ with $\iota_u(g)=ugu^{-1}$ for $u\in G$.

\begin{theorem}
Let $G$ be a separable TSI Lie group with $G_0=\R$. Define $\iota_\R:\R\to\R$ as $\iota_\R(t)=-t$. Then we have
\begin{enumerate}
\item[(1)] $E(G)\sim_BE(\R)\iff{\rm Inn}_G(G_0)=\{{\rm id}_\R\}$;
\item[(2)] $E(G)\sim_BE(\Z_2\ltimes\R)\iff{\rm Inn}_G(G_0)=\{{\rm id}_\R,\iota_\R\}$.
\end{enumerate}
\end{theorem}

Similar results of $G_0=\T$ in place of $\R$ are also valid.

In the end, the groups $\T_p$ of $p$-adic solenoids for all $p\ge 2$ are considered. Let $N(p)$ be the set of all prime divisors of $p$.

\begin{theorem}
Let $p,q\ge 2$ be natural numbers, then we have
\begin{enumerate}
\item[(1)] $E(\T_p)\le_BE(\T_q)\iff N(p)\supseteq N(q)$;
\item[(2)] $E(\R)<_BE(\T_p)<_BE(\T)$.
\end{enumerate}
\end{theorem}

The authors also generalized the above theorem to $P$-adic solenoids, where $P$ is a sequence of natural numbers $\ge 2$ (see~\cite[Theorem 3.2]{DZ}). Moreover, the Borel reducibility among $E(G)$'s between $E(\R)$ and $E(\T)$ are extremely complicated that the partial ordered set $P(\omega)/{\rm Fin}$ embeds into them (see~\cite[Theorem 3.6]{DZ}).
 

This article is organized as follows. In section 2, we recall some notions in descriptive set theory and Polish groups. In section 3, we define the equivalence relations induced by Polish groups, and prove Theorem 1.1 and some basic facts. In section 4, we prove theorems 1.2 and 1.3. In section 5, we define the notion of $\alpha$-l.m.-unbalanced groups and prove theorems 1.4 and 1.5. In section 6, we present many results on TSI Polish groups and prove theorems 1.6--1.10. Finally, in section 7, we give some further open questions and additional remarks.

\section{Preliminaries}

We say a topological space is {\it Polish} if it is separable and admits a compatible complete metric. Let $X,Y$ be two Polish spaces, $E$ and $F$ equivalence relations on $X$ and $Y$ respectively. We say a map $\theta:X\to Y$ is a {\it Borel reduction} of $E$ to $F$ if $\theta$ is Borel and for $x,y\in X$,
$$xEy\iff\theta(x)F\theta(y).$$
If such a Borel reduction exists, we say $E$ is {\it Borel reducible to} $F$, denoted by $E\le_BF$. We also denote by $E\sim_BF$ for $E\le_BF$ and $F\le_BE$; and denote by $E<_BF$ for $E\le_BF$ but $F\not\le_BE$. For more details on Borel reducibility and descriptive set theory, we refer to~\cite{gaobook,kanovei,kechris}.

The identity element of a group $G$ is denoted by $1_G$. We say a topological group is a {\it Polish group} if its underlying topology is Polish. Given a Polish group $G$ and a Polish space $X$, an action of $G$ on $X$, denoted by $G\curvearrowright X$, is a map $a:G\times X\to X$ satisfying that $a(1_G,x)=x$ and $a(gh,x)=a(g,a(h,x))$ for $g,h\in G$ and $x\in X$. The pair $(X,a)$ is called a {\it Polish~(Borel) $G$-space} if $a$ is a continuous~(Borel) map.
Throughout this article, we write $gx$ in place of $a(g,x)$ for the sake of brevity. For each $x\in X$, the {\it $G$-orbit} of $x$, denoted by $[x]$ or $Gx$, is the set $\{gx:g\in G\}$. A subset $A$ of $X$ is {\it $G$-invariant} if $Gx\subseteq A$ for each $x\in A$. The {\it orbit equivalence relation}, denoted by $E_G^X$, is defined as
$$xE_G^Xy\iff[x]=[y]\iff\exists g\in G\,(gx=y).$$

Becker and Kechris proved that, for any Borel $G$-space $(X,a)$, there is a Polish topology $\tau$ on $X$ such that $\tau$ generates the Borel structure of $X$ and $(X,a)$ is a Polish $G$-space~(cf.~\cite[Theorem 4.4.6]{gaobook}). This implies that $E^X_G=E^{(X,\tau)}_G$. Therefore, the concepts of Borel $G$-space and Polish $G$-space are equivalent under the perspective of Borel reducibility.

Note that any discrete countable group is Polish. The {\it infinite permutation group} $S_\infty$ is the group of all permutations of $\omega$ equipped with the pointwise convergence topology. A Polish group $G$ is {\it non-archimedean} if it has a neighborhood base of $1_G$ consisting of open subgroups. Becker and Kechris proved that a Polish group is non-archimedean iff it is topologically isomorphic to a closed subgroup of $S_\infty$ (cf.~\cite[Theorem 1.5.1]{BK} or~\cite[Theorem 2.4.1]{gaobook}). A metric $d$ on a group $G$ is {\it left-invariant} if $d(gh,gk)=d(h,k)$ for all $g,h,k\in G$. Similarly, we can define {\it right-invariant} metric. We say $d$ is {\it two-sided invariant} if it is both left and right-invariant. A compatible two-sided invariant metric on a Polish group is necessarily complete (cf.~\cite[Corollary 1.2.2]{BK}). We say a Polish group $G$ is CLI if it admits a left-invariant compatible complete metric; and say $G$ is TSI if it admits a two-sided invariant compatible metric.

Now we recall some benchmark equivalence relations in the research of Borel reducibility. The equivalence relation $E_0$ on $2^\omega$ is defined as
$$xE_0y\iff\exists m\,\forall n>m\,(x(n)=y(n)).$$
If $E$ is an equivalence relation on a Polish space $X$, then we define equivalence relations $E^\omega$ and $E^+$ on $X^\omega$ as
$$xE^\omega y\iff\forall n\,(x(n)Ey(n)),$$
 $$xE^+y\iff\forall n\,\exists m\,(x(n)Ey(m))\wedge\forall m\,\exists n\,(x(n)Ey(m)).
$$
We pay special attention to $E_0^\omega$ and ${\rm id}(\omega^\omega)^+$ (denoted by $=^+$ for brevity).

Let $G$ and $\Lambda$ be two groups, and $\phi$ a homomorphism from $\Lambda$ to ${\rm Aut}(G)$, the group of automorphisms on $G$. Recall that the {\it semi product} $\Lambda\ltimes_\phi G$ is the set $\Lambda\times G$ equipped with group operation as: for $(\lambda_1,g_1),(\lambda_2,g_2)\in\Lambda\times G$,
$$(\lambda_1,g_1)(\lambda_2,g_2)=(\lambda_1\lambda_2,g_1\phi(\lambda_1)(g_2)).$$
Note that, if $G$ is a Polish group, $\Lambda$ is a countable discrete group, and $\phi(\lambda)$ is a continuous automorphism on $G$ for each $\lambda$, then $\Lambda\ltimes_\phi G$ equipped with the product topology on $\Lambda\times G$ is also a Polish group and $G$ is topologically isomorphic to the open normal subgroup $\{1_\Lambda\}\times G$.
We also recall that the {\it wreath product} $\Lambda\wr G$ is the set $\Lambda\times G^\Lambda$ equipped with group operation as: for $(\lambda_1,\chi_1),(\lambda_2,\chi_2)\in\Lambda\times G^\Lambda$,
$$(\lambda_1,\chi_1)(\lambda_2,\chi_2)=(\lambda_1\lambda_2,\chi)$$
with $\chi(\lambda)=\chi_1(\lambda)\chi_2(\lambda_1^{-1}\lambda)$ for $\lambda\in\Lambda$. If $\Lambda$ is countable and $G$ is Polish, $\Lambda\wr G$ equipped the product topology of $\Lambda\times G^\Lambda$ is also a Polish group.

Given two sets $X,Y$ and a map $f:X\to Y$, we define a map $f^\omega:X^\omega\to Y^\omega$ as: for $x\in X^\omega$ and $n\in\omega$,
$$f^\omega(x)(n)=f(x(n)).$$

\section{Equivalence relations induced by Polish groups}

\begin{definition}
Let $G$ be a Polish group. We define an equivalence relation $E(G)$ on $G^\omega$ as
$$xE(G)y\iff\lim_nx(n)y(n)^{-1}\mbox{ converges in }G$$
for $x,y\in G^\omega$. We say $E(G)$ is the \textit{equivalence relation induced by $G$}. Moreover, we define a subgroup of $G^\omega$ as
$$c(G)=\{x\in G^\omega:\lim_nx(n)\mbox{ converges in }G\}.$$
Then we have
$$xE(G)y\iff xy^{-1}\in c(G)\iff c(G)x=c(G)y.$$
\end{definition}

It is easy to see that $E(G)$ is a ${\mathbf\Pi}^0_3$ equivalence relation on $G^\omega$.

Let $d$ be a compatible metric on $G$. We define the supremum metric as $d_u(x,y)=\sup\{d(x(n),y(n)):n\in\omega\}$ for $x,y\in G^\omega$. For any sequence $(x_k)$ in $G^\omega$ and any $x\in G^\omega$, we use $x_k\rightrightarrows x$ to denote $d_u(x_k,x)\to 0$, i.e., $(x_k)$ converges to $x$ uniformly (respect to $d$).

\begin{theorem}
Let $G$ be a Polish group, $d$ a complete compatible metric on $G$, and let $d_u$ be the supremum metric on $c(G)$. Then $(c(G),d_u)$ is also a Polish group.

Furthermore, let $d'$ be any compatible metric on $G$ and $d_u'$ the supremum metric of $d'$, then $(c(G),d_u)$ and $(c(G),d_u')$ have the same topology.
\end{theorem}

\begin{proof}
(1) First, we prove that $d_u$ is a Polish metric on $c(G)$. Let $(\gamma_m)$ be a $d_u$-Cauchy sequence in $c(G)$, and let $g_m=\lim_n\gamma_m(n)$. Then for each $n\in\omega$, $(\gamma_m(n))$ is a $d$-Cauchy sequence. So there exists an $x\in G^\omega$ such that $\lim_m\gamma_m(n)=x(n)$ for each $n\in\omega$. It is clear that $(g_m)$ is also a $d$-Cauchy sequence, so $g_m\to g$ for some $g\in G$. From the property of uniform convergence, we have
$$\lim_nx(n)=\lim_n\lim_m\gamma_m(n)=\lim_m\lim_n\gamma_m(n)=\lim_mg_m=g.$$
So $x\in c(G)$, and hence $d_u$ is complete.

(2) Secondly, we prove that $(c(G),d_u)$ is a topological group. Let $\gamma,\eta\in c(G)$ and let $(\gamma_m),(\eta_m)$ be two sequences in $c(G)$. Put $g=\lim_n\gamma(n)$, $h=\lim_n\eta(n)$, $g_m=\lim_n\gamma_m(n)$, and $h_m=\lim_n\eta_m(n)$ for each $m\in\omega$. Now suppose $d_u(\gamma_m,\gamma)\to 0$ and $d_u(\eta_m,\eta)\to 0$. We only need to show that $d_u(\gamma_m\eta_m^{-1},\gamma\eta^{-1})\to 0$.

For any $\varepsilon>0$, there exists a $\delta_1>0$ such that, for any $g',h'\in G$, if $d(g',g)<\delta_1$ and $d(h',h)<\delta_1$, then we have $d(g'(h')^{-1},gh^{-1})<\varepsilon/2$. For this $\delta_1>0$, there exists an $N_1\in\omega$ such that, for $n>N_1$, we have $d(\gamma(n),g)<\delta_1/2$ and $d(\eta(n),h)<\delta_1/2$. It follows that $d(\gamma(n)\eta(n)^{-1},gh^{-1})<\varepsilon/2$. There also exists an $M_1\in\omega$ such that, for $m>M_1$, we have $d_u(\gamma_m,\gamma)<\delta_1/2$ and $d_u(\eta_m,\eta)<\delta_1/2$. So
$$d(\gamma_m(n),g)\le d(\gamma_m(n),\gamma(n))+d(\gamma(n),g)\le d_u(\gamma_m,\gamma)+d(\gamma(n),g)<\delta_1,$$
and also $d(\eta_m(n),h)<\delta_1$. It follows that $d(\gamma_m(n)\eta_m(n)^{-1},gh^{-1})<\varepsilon/2$. Then we have
$$d(\gamma_m(n)\eta_m(n)^{-1},\gamma(n)\eta(n)^{-1})<\varepsilon\quad(\forall m>M_1,\forall n>N_1).$$
For this $\varepsilon>0$, there exist a $\delta_2>0$ such that, for any $g',h'\in G$ and any $n\le N_1$, if $d(g',\gamma(n))<\delta_2$ and $d(h',\eta(n))<\delta_2$, then $d(g'(h')^{-1},\gamma(n)\eta(n)^{-1})<\varepsilon$. Then there exists an $M_2\in\omega$ such that $d_u(\gamma_m,\gamma)<\delta_2$ and $d_u(\eta_m,\eta)<\delta_2$ for $m>M_2$.
Thus we have
$$d(\gamma_m(n)\eta_m(n)^{-1},\gamma(n)\eta(n)^{-1})<\varepsilon\quad(\forall m>M_2,\forall n\le N_1).$$
Therefore, $d_u(\gamma_m\eta_m^{-1},\gamma\eta^{-1})\le\varepsilon$ for all $m>\max\{M_1,M_2\}$. It follows that $d_u(\gamma_m\eta_m^{-1},\gamma\eta^{-1})\to 0$.

It is trivial that $(c(G),d_u)$ is separable, so it is a Polish group.

(3) Let $d'$ be an compatible metric on $G$ and $d'_u$ the supremum metric (respect to $d'$) on $c(G)$. Let $\gamma\in c(G)$ and let $(\gamma_m)$ be a sequence in $c(G)$, we only need to show that: $d_u(\gamma_m,\gamma)\to 0$ iff $d'_u(\gamma_m,\gamma)\to 0$. Let $\lim_n\gamma(n)=g$ and $\lim_n\gamma_m(n)=g_m$ for each $m\in\omega$. Suppose $d_u(\gamma_m,\gamma)\to 0$. Following similar arguments in the proof of $d_u(\gamma_m\eta_m^{-1},\gamma\eta^{-1})\to 0$ in (2), we can prove that $d'_u(\gamma_m,\gamma)\to 0$. And vice versa.
\end{proof}

Therefore, $E(G)$ is the orbit equivalence relation induced by $c(G)\curvearrowright G^\omega$, the action of left multiplication of $c(G)$ on $G^\omega$.

\begin{proposition}\label{CliTsiNonArchi}
Let $G$ be a Polish group. If $G$ is CLI (TSI, or non-archimedean), so is $c(G)$.
\end{proposition}

\begin{proof}
If $G$ is CLI, let $d$ be a left-invariant compatible complete metric on $G$. Then it is trivial to check that the supremum metric $d_u$ is also a left-invariant compatible complete metric on $c(G)$, so $c(G)$ is CLI too.

If $G$ is TSI, the same arguments show that $c(G)$ is TSI.

If $G$ is non-archimedean, let $(G_n)$ be a sequence of open subgroups of $G$ which forms a neighborhood base of $1_G$, and let $d$ be a compatible metric on $G$ and $d_u$ the supremum metric on $c(G)$. Find an $\varepsilon_n>0$ such that
$$\{g\in G:d(g,1_G)<\varepsilon_n\}\subseteq G_n.$$
Let $V_n=\{\gamma\in c(G):d_u(\gamma,1_{c(G)})<\varepsilon_n\}$. Then for any $\eta\in G_n^\omega\cap c(G)$, we have
$$\eta\in\eta V_n\subseteq G_n^\omega\cap c(G).$$
So $G_n^\omega\cap c(G)$ is an open subgroup of $c(G)$. It is clear that $(G_n^\omega\cap c(G))$ is a neighborhood base of $1_{c(G)}$, so $c(G)$ is non-archimedean.
\end{proof}

\begin{proposition}\label{subgroup}
Let $G,H$ be two Polish groups. If $G$ is topologically isomorphic to a closed subgroup of $H$, then $E(G)\le_BE(H)$.
\end{proposition}

\begin{proof}
It follows from the definitions of $E(G)$ and $E(H)$.
\end{proof}

\begin{theorem}\label{trivial}
\begin{enumerate}
\item[(1)] Let $G$ be a discrete countable group containing at least two elements, then $E(G)\sim_BE_0$.
\item[(2)] Let $G$ be an uncountable Polish group, then $E_0^\omega\le_BE(G)$.
\item[(3)] Let $G$ be a TSI non-archimedean uncountable Polish group, then $E(G)\sim_BE_0^\omega$.
\item[(4)] Let $G$ be a non-archimedean Polish group, then $E(G)\le_B=^+$.
\item[(5)] $E(S_\infty)\sim_B =^+$.
\end{enumerate}
\end{theorem}

\begin{proof}
(1)  Let $G$ be a discrete countable group containing at least two elements. For $x,y\in G^\omega$, we have
$$xE(G)y\iff xy^{-1}\in c(G)\iff\exists m\,\exists g\in G\,\forall n>m\,(x(n)y(n)^{-1}=g).$$
For all $m\in\omega$, we define
$$\begin{array}{ll}xF_my & \iff\exists g\in G\,\forall n>m\,(x(n)y(n)^{-1}=g)\cr & \iff\forall n>m\,(x(n)y(n)^{-1}=x(n+1)y(n+1)^{-1}).\end{array}$$
Thus each $F_m$ is a closed equivalence relation on $G^\omega$, and hence is smooth (cf. \cite[Proposition 5.4.7]{gaobook}). Note that $E(G)=\bigcup_mF_m$ and $F_m\subseteq F_{m+1}$ for each $m\in\omega$. So $E(G)$ is hypersmmoth (cf. \cite[Definition 8.1.1]{gaobook}). It is clear that each orbit of $E(G)$ is countable. By Dougherty-Jackson-Kechris' theorem (cf. \cite[Theorem 8.1.5]{gaobook}), we have $E(G)\le_BE_0$.

Fix a $g_0\in G$ with $g_0\ne 1_G$. For $a\in 2^\omega$, we define $\theta(a)\in G^\omega$ as
$$\theta(a)(n)=\left\{\begin{array}{ll}g_0, & n=2k,a(k)=1,\cr 1_G, &\mbox{otherwise.}\end{array}\right.$$
Then $\theta$ witnesses that $E_0\le_BE(G)$. Therefore, we have $E(G)\sim_BE_0$.

(2) Let $G$ be uncountable, then we can find a sequence $(g_i)$ in $G$ such that $\lim_i g_i=1_G$ with all $g_i\ne 1_G$ for $i\in\omega$. Fix a bijection $\langle\cdot,\cdot\rangle$ from $\omega\times\omega$ to $\omega$. For $\alpha\in(2^\omega)^\omega$, we define $\vartheta(\alpha)\in G^\omega$ as
$$\vartheta(\alpha)(\langle i,j\rangle)=\left\{\begin{array}{ll}g_i, & j=2k,\alpha(i)(k)=1,\cr 1_G, &\mbox{otherwise.}\end{array}\right.$$
Then $\vartheta$ witnesses that $E_0^\omega\le_BE(G)$.

(3) Let $G$ be TSI non-archimedean uncountable. By (2), we only need to show $E(G)\le_BE_0^\omega$. From \cite[Theorem 1.1]{GX}, $G$ is isomorphic to a closed subgroup of a product $\prod_k\Gamma_k$, where each $\Gamma_k$ is a discrete countable group. So, by Proposition~\ref{subgroup}, it suffices to prove $E(\prod_k\Gamma_k)\le_BE_0^\omega$. For $x,y\in(\prod_k\Gamma_k)^\omega$, we have
$$\begin{array}{ll}xE(\prod_k\Gamma_k)y &\iff xy^{-1}\in c(\prod_k\Gamma_k)\cr
&\iff\forall k\,(x(\cdot)(k)y(\cdot)(k)^{-1}\in c(\Gamma_k))\cr
&\iff\forall k\,(x(\cdot)(k)E(\Gamma_k)y(\cdot)(k)).\end{array}$$
Then $E(\prod_k\Gamma_k)\le_BE_0^\omega$ follows from (1).

(4) From \cite[Theorem 2.4.1]{gaobook} and Proposition~\ref{subgroup}, it suffices to prove $E(S_\infty)\le_B=^+$. For $x,y\in (S_\infty)^\omega$, we have
$$\begin{array}{ll}xE(S_\infty)y &\iff xy^{-1}\in c(S_\infty)\cr
&\iff\exists g\in S_\infty\,\forall k\,\exists m\,\forall n>m\,((x(n)y(n)^{-1})(k)=g(k))\cr
&\iff\exists g\in S_\infty\,\forall k\,\exists m\,\forall n>m\,(x(n)^{-1}(g(k))=y(n)^{-1}(k)).\end{array}$$
We define $\theta'(x)\in(\omega^\omega)^\omega$ as $\theta'(x)(k)(n)=x(n)^{-1}(k)$ for $n,k\in\omega$. Then
$$xE(S_\infty)y\iff\exists g\in S_\infty\,\forall k\,(\theta'(x)(g(k))E_0\theta'(y)(k)).$$
By \cite[Exercise 8.3.4]{gaobook}, this implies that $E(S_\infty)\le_BE_0^+$. It is well known that $E_0^+\le_B=^+$ (cf. \cite[Theoerm 7.4.10 and Exercise 8.3.3]{gaobook}), so we have $E(S_\infty)\le_B=^+$.

(5) By~\cite[Lemma 10.3.4]{gaobook} and (4), we only need to show that $=^+\upharpoonright Y\le_BE(S_\infty)$,
where $Y=\{\alpha\in (2^\omega)^\omega:\forall n,m\,(n\ne m\Rightarrow\alpha(n)\ne\alpha(m))\}$.
Let $(p_i)$ be the strictly increasing enumeration of all prime numbers. Define a map $\varphi:Y\rightarrow(\omega^\omega)^\omega$ with
$$\varphi(\alpha)(n)(k)=\left\{\begin{array}{ll} p_0^{\alpha(n)(0)}\cdots p_k^{\alpha(n)(k)}, & k\ge n,\cr 0, &\mbox{otherwise},\end{array}\right.$$
for $\alpha\in Y$ and $n,k\in\omega$. It is clear that
$$\varphi(\alpha)(n)E_0\varphi(\beta)(m)\iff\alpha(n)=\beta(m).$$

Define a map $\rho:(\omega^\omega)^\omega\rightarrow(\omega^\omega)^\omega$ by induction on $n$ as:
$$\rho(x)(k)(n)=\left\{\begin{array}{ll} x(n)(k), & x(n)(k)\notin\{\rho(x)(k)(i):i<n\},\cr \min{\omega\setminus\{\rho(x)(k)(i):i<n\}}, &\mbox{otherwise},\end{array}\right.$$
for $x\in(\omega^\omega)^\omega$ and $n,k\in\omega$. Note that $\varphi(\alpha)(n)(k)=0$ for all $n>k$, so $\rho(\varphi(\alpha))(k)\in S_\infty$ for all $\alpha\in Y$ and $k\in\omega$.
Furthermore, for any $n\in\omega$, $\alpha(0),\dots,\alpha(n)$ are pairwise distinct elements in $2^\omega$, so for large enough $k$,
we have that $\varphi(\alpha)(0)(k),\dots,\varphi(\alpha)(n)(k)$ are pairwise distinct natural numbers. It follows that
$$\forall n\,\forall^\infty k\,(\rho(\varphi(\alpha))(k)(n)=\varphi(\alpha)(n)(k)).$$

Now we define a map $\vartheta':Y\rightarrow(S_\infty)^\omega$ as
$$\vartheta'(\alpha)(k)=\rho(\varphi(\alpha))(k)^{-1}\in S_\infty$$
for $x\in Y$ and $k\in\omega$. Then for any $\alpha,\beta\in Y$, we have
$$\begin{aligned}
&\vartheta'(\alpha)E(S_\infty)\vartheta'(\beta)\cr
\iff &\exists g\in S_\infty\,\forall n\,\forall^\infty k\,((\rho(\varphi(\alpha))(k)^{-1}\rho(\varphi(\beta))(k))(n)=g(n))\cr
\iff &\exists g\in S_\infty\,\forall n\,\forall^\infty k\,(\rho(\varphi(\alpha))(k)(g(n))=\rho(\varphi(\beta))(k)(n))\cr
\iff &\exists g\in S_\infty\,\forall n\,\forall^\infty k\,(\varphi(\alpha)(g(n))(k)=\varphi(\beta)(n)(k))\cr
\iff &\exists g\in S_\infty\,\forall n\,(\varphi(\alpha)(g(n))E_0\varphi(\beta)(n))\cr
\iff &\exists g\in S_\infty\,\forall n\,(\alpha(g(n))=\beta(n))\cr
\iff &\alpha=^+\beta.
\end{aligned}
$$
Thus $=^+\upharpoonright Y\le_B E(S_\infty)$. So $E(S_\infty)\sim_B=^+$.
\end{proof}

Furthermore, the authors~\cite{DZtsi} get the following stronger result.

\begin{thm}[{\cite[Theorem 1.3]{DZtsi}}]
Let $G$ be a Polish group. Then the following are equivalent:
\begin{enumerate}
  \item $G$ is TSI non-archimedean;
  \item $E(G)\leq_B E_0^\omega$; and
  \item $E(G)\leq_B \R^\omega/c_0$.
\end{enumerate}
In particular, $E(G)\sim_B E_0^\omega$ iff $G$ is TSI uncountable non-archimedean.
\end{thm}

\begin{corollary}\label{timesE_0}
Let $\Lambda$ be a discrete countable group, $G$ a Polish group containing an element $h_0$ with $1_G\notin\overline{\{gh_0g^{-1}:g\in G\}}$. Then $E(G\times\Lambda)\sim_BE(G)$.
\end{corollary}

\begin{proof}
Proposition~\ref{subgroup} gives $E(G)\le_BE(G\times\Lambda)$. It is trivial that $E(G\times\Lambda)\sim_BE(G)\times E(\Lambda)$. By Theorem~\ref{trivial}(1), $E(\Lambda)\le_BE_0$, so we only need to show that $E(G)\times E_0\le_BE(G)$. For $x\in G^\omega$ and $a\in 2^\omega$, define $\theta(x,a)\in G^\omega$ as
$$\theta(x,a)(n)=\left\{\begin{array}{ll}x(k)h_0, & n=2k,a(k)=1,\cr x(k), & n=2k,a(k)=0\mbox{ or }n=2k+1.\end{array}\right.$$
It is trivial to check that $\theta$ is a Borel reduction of $E(G)\times E_0$ to $E(G)$.
\end{proof}

\section{Non-CLI Polish groups and Non Non-archimedean Polish groups}

\subsection{Non-CLI Polish groups}
Lupini and Panagiotopoulos introduced a notion of right $\iota$-embedability for proving that an orbit equivalence relation is not Borel reducible to any $E_H^Y$, where $H$ is a CLI Polish group.

\begin{definition}[Becker \cite{becker}, Lupini--Panagiotopoulos \cite{LP}]
Let $G$ be a Polish group, $(g_n)$ a sequence in $G$. We say that $(g_n)$ is \textit{left Cauchy} if $\lim_{m,n}g_m^{-1}g_n=1_G$; and say that $(g_n)$ is \textit{right Cauchy} if $\lim_{m,n}g_mg_n^{-1}=1_G$.

Let $X$ be a Polish $G$-space and $x,y\in X$. We say that $x$ is \textit{left $\iota$-embeddable} into $y$, denoted by $x\stackrel{l}{\to}y$, if there exists a left Cauchy sequence $(g_n)$ in $G$ such that $g_nx\to y$; and say that $x$ is \textit{right $\iota$-embeddable} into $y$, denoted by $x\stackrel{r}{\to}y$, if there exists a right Cauchy sequence $(h_n)$ in $G$ such that $h_ny\to x$.
\end{definition}

Clearly, $(g_n)$ is right Cauchy iff $(g_n^{-1})$ is left Cauchy. Let $d$ be a left-invariant compatible metric on $G$, $(g_n)$ a sequence in $G$. It is worth noting that, $(g_n)$ is left Cauchy iff it is $d$-Cauchy. It was proved by Becker that the relation of left $\iota$-embeddability is a preorder, and is invariant (cf. \cite{becker}). Similar arguments show that the relation of right $\iota$-embeddability is a preorder, and is invariant too. So we can write $[x]\stackrel{l}{\to}[y]$ and $[x]\stackrel{r}{\to}[y]$ in place of $x\stackrel{l}{\to}y$ and $x\stackrel{r}{\to}y$ respectively.

\begin{theorem}
Let $X$ be a Polish $G$-space and $x,y\in X$. Then $x\stackrel{l}{\to}y$ implies $x\stackrel{r}{\to}y$.
\end{theorem}

\begin{proof}
Suppose $x\stackrel{l}{\to}y$. Then there exists a left Cauchy sequence $(g_n)$ in $G$ such that $g_nx\to y$. Let $d$ be a compatible metric on $G$ and $\rho$ a compatible metric on $X$.

By the continuity of the action of $G$ on $X$, for any $k\in\omega$, there exists a $\delta>0$ such that, for any $g\in G$, if $d(g,1_G)<\delta$, then $\rho(gx,x)<2^{-(k+1)}$.
Since $(g_n)$ is left Cauchy, for this $\delta>0$, there exists an $m_k\in\omega$ such that, for any $n>m_k$, we have $d(g_{m_k}^{-1}g_n,1_G)<\delta$, and hence $\rho(g_{m_k}^{-1}g_nx,x)<2^{-(k+1)}$. By $g_nx\to y$, we have $g_{m_k}^{-1}g_nx\to g_{m_k}^{-1}y$. There exists an $n>m_k$ such that $\rho(g_{m_k}^{-1}g_nx,g_{m_k}^{-1}y)<2^{-(k+1)}$. So we have $\rho(g_{m_k}^{-1}y,x)<2^{-k}$. Without loss of generality, we can assume that $(m_k)$ is strictly increasing. Let $h_k=g_{m_k}^{-1}$ for each $k$. Then $(h_k)$ witnesses that $x\stackrel{r}{\to}y$.
\end{proof}

\begin{theorem}\label{CLI}
Let $G$ be a non-CLI Polish group and $H$ a CLI Polish group. Then for any Polish $H$-space $Y$, we have $E(G)\not\le_BE_H^Y$.
\end{theorem}

\begin{proof}
Let $d$ be a left-invariant compatible metric on $G$, and let $d_u$ be the supremum metric on $c(G)$. Since $G$ is not CLI, $d$ is not complete. So there exists a $d$-Cauchy sequence $(g_n)$ in $G$ such that $(g_n)$ diverges. For $m\in\omega$, we define $\gamma_m\in c(G)$ as:
$$\gamma_m(n)=\left\{\begin{array}{ll}g_n, & n<m,\cr g_m, & n\ge m.\end{array}\right.$$
For $m<k$, we have $d_u(\gamma_m^{-1}\gamma_k,1_{c(G)})=\max\{d(g_m,g_n):m\le n\le k\}$. So $(\gamma_m)$ is left Cauchy, and hence $(\gamma_m^{-1})$ is right Cauchy.

We define $z\in G^\omega$ as $z(n)=g_n$ for each $n\in\omega$. Then $z\notin c(G)$. For any $c(G)$-invariant dense $G_\delta$ subset $C$ of $G^\omega$, since $z^{-1}C$ is also dense $G_\delta$, $C\cap z^{-1}C\ne\emptyset$. Let $x\in C\cap z^{-1}C$ and $y=zx$, then $x,y\in C$ and $[x]\ne[y]$. Note that $\gamma_m^{-1}y\to z^{-1}y=x$, so $x\stackrel{r}{\to} y$. In the end, \cite[Lemma 2.4 and Theorem 2.9]{LP} gives $E(G)\not\le_BE_H^Y$.
\end{proof}

\begin{corollary}
Let $H$ be a CLI Polish group, then $E(S_\infty)\not\le_BE(H)$.
\end{corollary}

\begin{proof}
It well known that $S_\infty$ is not CLI (cf.~\cite[Example 2.2.7]{gaobook}). And by Proposition~\ref{CliTsiNonArchi}, $c(H)$ is CLI.
\end{proof}

\subsection{Non Non-archimedean Polish groups}

Hjorth founded the theory of turbulence for proving that an orbit equivalence relation is not Borel reducible to $E_{S_\infty}^Y$ for any Polish $S_\infty$-space $Y$. We omit the definition of turbulent actions since it is very complexity and is not used in this article, one can find it in~\cite[\S 10]{gaobook}.

\begin{definition}[Farah~\cite{farah}]
If $(X_n,d_n),\,n\in\omega$ is a sequence of finite metric spaces, let $D=D(\langle X_n,d_n\rangle)$ be the equivalence relation on $\prod_nX_n$ defined as
$$xDy\iff\lim_nd_n(x(n),y(n))=0.$$
The equivalence relations of this from are called {\it $c_0$-equalities}.
\end{definition}

By~\cite[Lemma 3.4]{farah}, every $c_0$-equality is an orbit equivalence relation induced by a continuous Polish group action on $\prod_nX_n$.

For a finite metric space $(X,d)$ and $\varepsilon\in\R_+$, we define:
$$\delta(\varepsilon,X)=\min\{\delta:\exists x_0,x_1,\cdots,x_n\,(d(x_0,x_n)\ge\varepsilon\wedge\forall i<n\,(d(x_i,x_{i+1})<\delta))\}.$$

\begin{lemma}
Let $G$ be a Polish group. If $G$ is not non-archimedean, then there exists a turbulent $c_0$-equality $D$ such that $D\le_BE(G)$.
\end{lemma}

\begin{proof}
Let $d$ be a left-invariant compatible metric on $G$. Since $G$ is not non-archimedean, there exists an $\varepsilon_0>0$ such that no open subgroup of $G$ is contained in $V=\{g\in G:d(1_G,g)<\varepsilon_0\}$. For each $n\in\omega$, we define $V_n=\{g\in G:d(1_G,g)<2^{-n}\}$. Then $V_n^{-1}=V_n$ and $G_n=\bigcup_mV_n^m$ is an open subgroup of $G$. Thus $G_n\not\subseteq V$ for each $n\in\omega$, and hence there is an integer $m_n>0$ such that $V_n^{m_n}\not\subseteq V$. So we can find $1_G=g_{n,0},g_{n,1},\cdots,g_{n,m_n}\in G$ such that $g_{n,i}^{-1}g_{n,i+1}\in V_n$ for each $i<m_n$ and $g_{n,m_n}\notin V$.

Let $X_n=\{g_{n,i}:0\le i\le m_n\}$ for each $n\in\omega$. Then $(X_n,d)$ is a finite metric space. We define $\theta:\prod_nX_n\to G^\omega$ as, for $x\in\prod_nX_n$ and each $k\in\omega$,
$$\theta(x)(k)=\left\{\begin{array}{ll}x(n)^{-1}, & k=2n,\cr 1_G, & k=2n+1.\end{array}\right.$$
It is trivial to check that $\theta$ is a Borel reduction of $D(X_n,d)$ to $E(G)$.

Note that $\delta(\varepsilon_0,X_n)\le 2^{-n}\to 0$ as $n\to\infty$. By~\cite[Theorem 3.7(a)]{farah}, there exists a turbulent $c_0$-equality $D$ such that $D\le_BD(X_n,d)\le_BE(G)$.
\end{proof}

\begin{theorem}
Let $G$ and $H$ be two Polish groups. If $H$ is non-archimedean but $G$ is not, then $E(G)\not\le_BE_H^Y$ for any Polish $H$-space $Y$.

In particular, we have $E(G)\not\le_BE(S_\infty)$.
\end{theorem}

\begin{proof}
Since $H$ is topologically isomorphic to a closed subgroup of $S_\infty$, by~\cite[Theorem 3.5.2]{gaobook}, we only need to consider the case of $H=S_\infty$. Then the result follows from the preceding lemma and the Hjorth turbulence theorem (cf.~\cite[Corollary 10.4.3]{gaobook}).

By Proposition~\ref{CliTsiNonArchi}, $c(S_\infty)$ is non-archimedean, so $E(G)\not\le_BE(S_\infty)$.
\end{proof}

\begin{corollary}
Let $G$ be a Polish group, then $E(G)\le_B=^+$ iff $G$ is non-archimedean.
\end{corollary}

\begin{proof}
It follows from Theorem~\ref{trivial}(4) and the preceding theorem.
\end{proof}

\section{Non-TSI Polish groups and $\alpha$-unbalanced relations}

\subsection{Unbalanced relations}

\begin{definition}[Allison--Panagiotopoulos {\cite[Definition 1.1]{AP}}]
Let $X$ be a Polish $G$-space and let $x,y\in X$. We write $x\leftrightsquigarrow y$ provided that, for any open $V\ni 1_G$ and open set $U$ of $X$ with $U\cap([x]\cup[y])\ne\emptyset$, there exist $g^x,g^y\in G$ with $g^xx,g^yy\in U$, so that
$$g^yy\in\overline{V(g^xx)}\mbox{ and }g^xx\in\overline{V(g^yy)}.$$

It is clear that $\leftrightsquigarrow$ is symmetric, and for any $g,h\in G$, we have that $x\leftrightsquigarrow y\iff gx\leftrightsquigarrow hy$. So we can write $[x]\leftrightsquigarrow[y]$ whenever $x\leftrightsquigarrow y$.
\end{definition}

\begin{lemma}\label{unbalanced}
Let $G$ be a Polish group and $c(G)\curvearrowright G^\omega$ the action of left multiplication. For $x,y\in G^\omega$ we have
$$x\leftrightsquigarrow y\iff\exists(\gamma_m),(\eta_m)\in c(G)^\omega\,(\gamma_mxy^{-1}\eta_m\rightrightarrows 1_{G^\omega}).$$
\end{lemma}

\begin{proof}
Let $d$ be a compatible metric on $G$ and $d_u$ the supremum metric on $G^\omega$. Note that $1_{c(G)}(n)=1_{G^\omega}(n)=1_G$ for each $n\in\omega$. Let $V_m=\{\gamma\in c(G):d_u(\gamma,1_{c(G)})<2^{-m}\}$ for $m\in\omega$. Then $(V_m)$ is a neighborhood basis of $1_{c(G)}$. Note that $V_m$ is dense in $\{x\in G^\omega:d_u(x,1_{G^\omega})<2^{-m}\}$ under the topology of $G^\omega$. For $x,y\in G^\omega$, we have $x\in\overline{V_my}\iff xy^{-1}\in\overline{V_m}$, so
$$d_u(xy^{-1},1_{G^\omega})<2^{-m}\Rightarrow x\in\overline{V_my}\Rightarrow d_u(xy^{-1},1_{G^\omega})\le 2^{-m}.$$

($\Rightarrow$). Let $x\leftrightsquigarrow y$. Then for each $m\in\omega$, there exist $\gamma^x_m,\gamma^y_m\in c(G)$ such that $\gamma^x_mx\in\overline{V_m(\gamma^y_my)}$. So $d_u(\gamma^x_mxy^{-1}(\gamma^y_m)^{-1},1_{G^\omega})\le 2^{-m}$. Let $\gamma_m=\gamma^x_m$ and $\eta_m=(\gamma^y_m)^{-1}$, then $\gamma_mxy^{-1}\eta_m\rightrightarrows 1_{G^\omega}$.

($\Leftarrow$). Let $(\gamma_m),(\eta_m)\in c(G)^\omega$ with $\gamma_mxy^{-1}\eta_m\rightrightarrows 1_{G^\omega}$. It is easy to check that $\eta_m^{-1}yx^{-1}\gamma_m^{-1}\rightrightarrows 1_{G^\omega}$. For any open $V\ni 1_{c(G)}$ and open set $U$ of $G^\omega$ with $U\cap([x]\cup[y])\ne\emptyset$, there exists $k\in\omega$ with $V_k\subseteq V$. Without loss of generality, we may assume that $U\cap[x]\ne\emptyset$ and $U=U_0\times\cdots\times U_l\times G^\omega$ with $U_0,\cdots,U_l$ open in $G$. So there exists $\gamma\in c(G)$ with $\gamma x\in U$, i.e., $\gamma(n)x(n)\in U_n$ for each $n\le l$.

Fix a large enough $m$ such that
$$d_u(\gamma_mxy^{-1}\eta_m,1_{G^\omega})<2^{-k},\quad d_u(\eta_m^{-1}yx^{-1}\gamma_m^{-1},1_{G^\omega})<2^{-k}.$$
Now we define
$$\gamma^x(n)=\left\{\begin{array}{ll}\gamma(n), & n\le l,\cr \gamma_m(n), & n>l,\end{array}\right.
\quad\gamma^y(n)=\left\{\begin{array}{ll}\gamma(n)x(n)y(n)^{-1}, & n\le l,\cr \eta_m(n)^{-1}, & n>l.\end{array}\right.$$
Then we have $\gamma^x,\gamma^y\in c(G)$ and $\gamma^xx,\gamma^yy\in U$. Note that
$$\gamma^x(n)x(n)y(n)^{-1}(\gamma^y(n))^{-1}=\left\{\begin{array}{ll}1_G, & n\le l,\cr \gamma_m(n)x(n)y(n)^{-1}\eta_m(n), & n>l.\end{array}\right.$$
It is easy to see that
$$d_u(\gamma^xxy^{-1}(\gamma^y)^{-1},1_{G^\omega})\le d_u(\gamma_mxy^{-1}\eta_m,1_{G^\omega})<2^{-k},$$
$$d_u(\gamma^yyx^{-1}(\gamma^x)^{-1},1_{G^\omega})\le d_u(\eta_m^{-1}yx^{-1}\gamma_m^{-1},1_{G^\omega})<2^{-k}.$$
It gives that $\gamma^xx\in\overline{V_k(\gamma^yy)}\subseteq\overline{V(\gamma^yy)}$ and $\gamma^yy\in\overline{V_k(\gamma^xx)}\subseteq\overline{V(\gamma^xx)}$. Thus we have $x\leftrightsquigarrow y$.
\end{proof}

\begin{definition}
Let $G$ be a Polish group and $c(G)\curvearrowright G^\omega$ the action of left multiplication. We say $G$ is {\it l.m.-unbalanced} if there exist $x,y\in G^\omega$ such that $x\leftrightsquigarrow y$ and $[x]\ne[y]$, where l.m. stands for left multiplication.
\end{definition}

\begin{theorem}\label{TSI}
Any TSI Polish group is not l.m.-unbalanced.
\end{theorem}

\begin{proof}
Let $G$ be a Polish group and $c(G)\curvearrowright G^\omega$ the action of left multiplication. Assume that $x\leftrightsquigarrow y$ for some $x,y\in G^\omega$. From Lemma~\ref{unbalanced}, there exist $(\gamma_m),(\eta_m)\in c(G)^\omega$ with $\gamma_mxy^{-1}\eta_m\rightrightarrows 1_{G^\omega}$. Let $d$ be a two-sided invariant compatible metric on $G$ and $d_u$ the supremum metric on $G^\omega$. Then $d_u$ is also two-sided invariant, so
$$d_u(\eta_m\gamma_mxy^{-1},1_{G^\omega})=d_u(\gamma_mxy^{-1}\eta_m,1_{G^\omega})\to 0,$$
thus we have $\eta_m\gamma_mxy^{-1}\rightrightarrows 1_{G^\omega}$. It follows that $\eta_m\gamma_m\rightrightarrows yx^{-1}$. Since $c(G)$ is closed under uniform convergence, $yx^{-1}\in c(G)$. So $[x]=[y]$.
\end{proof}

As an application of the notion of l.m.-unbalanced groups, we present the following theorem. Since it is a special case of Theorem~\ref{alpha-unbalanced}, we omit the proof at this moment.

\begin{theorem}\label{unbalanced to not}
Let $G,H$ be Polish groups. If $G$ is l.m.-unbalanced but $H$ is not, then $E(G)\not\le_BE(H)$.
\end{theorem}

\begin{corollary}
Let $G$ be a Polish group. If $G$ is l.m.-unbalanced, then $E_0^\omega<_BE(G)$.
\end{corollary}

\begin{proof}
It follows from Theorems~\ref{trivial}(2) and Theorem~\ref{unbalanced to not}.
\end{proof}

\subsection{$\alpha$-unbalanced relations}

\begin{definition}[Allison--Panagiotopoulos {\cite[Definition 7.1]{AP}}]
Let $X$ be a Polish $G$-space, $V$ an open neighborhood of $1_G$, and let $\alpha<\omega_1$. We define $\leftrightsquigarrow_V^\alpha$ by induction. We say
\begin{enumerate}
\item[(1)] $x\leftrightsquigarrow_V^0y$, if $y\in\overline{Vx}$ and $x\in\overline{Vy}$;
\item[(2)] $x\leftrightsquigarrow_V^\alpha y$ for $\alpha>0$, if for any open neighborhood $W\ni 1_G$ and any open neighborhood $U\subseteq X$ of $x$ or $y$, there exist $g^x,g^y\in V$ with $g^xx,g^yy\in U$, so that $g^xx\leftrightsquigarrow_W^\beta g^yy$ for all $\beta<\alpha$.
\end{enumerate}
\end{definition}

\begin{remark}
This definition appeared in an old version of~\cite{AP}. It is slightly different with which in new version of~\cite{AP} that the requirement of $U$ has been modified to $U\cap([x]\cup[y])\ne\emptyset$. It should be noted that the definitions of these two versions are not equivalent for $V\ne G$ or $\alpha>1$.

It is trivial that, if $x\leftrightsquigarrow_V^\alpha y$, then we have $x\leftrightsquigarrow_V^{\alpha'} y$ for $0<\alpha'\le\alpha$. Note that $g^xx\in U\cap Vx$ for any open neighborhood $U$ of $y$, so $U\cap Vx\ne\emptyset$, and hence $y\in\overline{Vx}$. Similarly, $x\in\overline{Vy}$. Therefore, $x\leftrightsquigarrow_V^\alpha y$ implies $x\leftrightsquigarrow_V^0 y$ too.
\end{remark}

\begin{lemma}~\label{gVh}
Let $V$ and $V'$ be two open neighborhoods of $1_G$ and $g,h\in G$ with $(gVg^{-1}\cup gVh^{-1}\cup hVh^{-1}\cup hVg^{-1})\subseteq V'$, and let $\alpha<\omega_1$, $x,y\in X$. If $x\leftrightsquigarrow_V^\alpha y$, then $gx\leftrightsquigarrow_{V'}^\alpha hy$. In particular, we have $gx\leftrightsquigarrow_{gVg^{-1}}^\alpha gy$.
\end{lemma}

\begin{proof}
We prove by induction on $\alpha$. If $\alpha=0$, since $x\leftrightsquigarrow_V^0y$ means $x\in\overline{Vy}$ and $y\in\overline{Vx}$, we have $gx\in\overline{gVy}=\overline{(gVh^{-1})(hy)}\subseteq\overline{V'(hy)}$ and $hy\in\overline{hVx}=\overline{(hVg^{-1})(gx)}\subseteq\overline{V'(gx)}$, so $gx\leftrightsquigarrow_{V'}^0hy$.

Assume that the result holds for all $\beta<\alpha$. If $x\leftrightsquigarrow_V^\alpha y$, then for any open $W\ni 1_G$ and any open set $U\subseteq X$ such that $U$ contains $gx$ or $hy$. Without loss of generality, suppose $gx\in U$, i.e., $x\in g^{-1}U$. Then there exist $g^x,g^y\in V$ with $g^xx,g^yy\in g^{-1}U$, so that $g^xx\leftrightsquigarrow_{g^{-1}Wg}^\beta g^yy$ for all $\beta<\alpha$. By the inductive hypothesis, $(gg^xg^{-1})(gx)=gg^xx\leftrightsquigarrow_W^\beta gg^yy=(gg^yh^{-1})(hy)$ for all $\beta<\alpha$. It is trivial to check that $gg^xg^{-1}$ and $gg^yh^{-1}$ witnesses that $gx\leftrightsquigarrow_{V'}^\alpha hy$.
\end{proof}

\begin{lemma}\label{1-unbalanced}
For all $\alpha<\omega_1$, the relation $\leftrightsquigarrow_G^\alpha$ is $G$-invariant. Moreover, $x\leftrightsquigarrow y$ iff $x\leftrightsquigarrow_G^1y$ for $x,y\in X$.
\end{lemma}

\begin{proof}
For any $x,y\in X$ and $g,h\in G$, by Lemma~\ref{gVh}, $x\leftrightsquigarrow_G^\alpha y$ implies $gx\leftrightsquigarrow_G^\alpha hy$. So $\leftrightsquigarrow_G^\alpha$ is $G$-invariant.

By their definitions, it is trivial that $x\leftrightsquigarrow y$ implies $x\leftrightsquigarrow_G^1y$.

On the other hand, if $x\leftrightsquigarrow_G^1y$, we will show $x\leftrightsquigarrow y$. For any open $V\ni 1_G$ and any open set $U\subseteq X$ with $U\cap([x]\cup[y])\ne\emptyset$, without loss of generality, we can assume that $U\cap[x]\ne\emptyset$. There exists $g\in G$ with $gx\in U$, i.e., $x\in g^{-1}U$. So there are $g^x,g^y\in G$ such that $g^xx,g^yy\in g^{-1}U$, $g^xx\in\overline{(g^{-1}Vg)g^yy}$, and $g^yy\in\overline{(g^{-1}Vg)g^xx}$. Then $gg^x$ and $gg^y$ witnesses that $x\leftrightsquigarrow y$.
\end{proof}

\begin{lemma}\label{V^xV^y}
Let $0<\alpha<\omega_1$, $V$ an open neighborhood of $1_G$. If $x\leftrightsquigarrow_V^\alpha y$, then for any open set $W\ni 1_G$ and any open set $U$ of $X$ such that $U$ contains $x$ or $y$, there exist nonempty open sets $V^x,V^y\subseteq V$ such that, whenever $g^x\in V^x$ and $g^y\in V^y$, we have $g^xx,g^yy\in U$ and $g^xx\leftrightsquigarrow_W^\beta g^yy$ for all $\beta<\alpha$.
\end{lemma}

\begin{proof}
Given an open set $W\ni 1_G$ and an open set $U$ of $X$ with $x\in U$ or $y\in U$. We can find an open $W_0\ni 1_G$ with $W_0^{-1}=W_0$ and $W_0^3\subseteq W$.
Then there exist $h^x,h^y\in V$ such that $h^xx,h^yy\in U$ and $h^xx\leftrightsquigarrow_{W_0}^\beta h^yy$ for all $\beta<\alpha$. We can find an open $V_0\ni 1_G$ such that $(V_0h^xx\cup V_0h^yy)\subseteq U$. By Lemma~\ref{gVh}, $w^xh^xx\leftrightsquigarrow_W^\beta w^yh^yy$ holds for all $w^x,w^y\in W_0$. Then we can put $V^x=W_0h^x\cap V_0h^x\cap V$ and $V^y=W_0h^y\cap V_0h^y\cap V$ as desired.
\end{proof}

\begin{theorem}\label{homomorphism}
Let $G,H$ be two Polish groups, $X$ a Polish $G$-space, and $Y$ a Polish $H$-space, and let $f:X\to Y$ be a Baire measurable homomorphism, i.e., $xE_G^Xy$ implies $f(x)E_H^Yf(y)$ for any $x,y\in X$. Then there exists a dense $G_\delta$ subset $C\subseteq X$ such that $x\leftrightsquigarrow_G^\alpha y$ implies $f(x)\leftrightsquigarrow_H^\alpha f(y)$ for any $x,y\in C$ and $0<\alpha<\omega_1$.
\end{theorem}

\begin{proof}
Let $C$ be the dense $G_\delta$ subset of $X$ in \cite[Lemma 2.5]{LP}.

{\bf Claim}. Let $0<\alpha<\omega_1$, $x_0,y_0\in C$, and let $V\ni 1_G$ open in $G$, $W\ni 1_H$ open in $H$, and $U$ open in $X$ with $x_0,y_0\in U$. Suppose
\begin{enumerate}
\item[(1)] $\forall x\in(U\cap C)\,\forall^*g\in V\,(f(gx)\in Wf(x))$;
\item[(2)] $x_0\leftrightsquigarrow_V^\alpha y_0$.
\end{enumerate}
Then we have $f(x_0)\leftrightsquigarrow_W^\alpha f(y_0)$.

{\it Proof of Claim}. By \cite[Lemma 2.5]{LP}, the following conditions hold:
\begin{enumerate}
\item[(i)] $f\upharpoonright C$ is continuous;
\item[(ii)] $\forall x\in C\,\forall^*g\in G\,(gx\in C)$;
\item[(iii)] for any open $W'\ni 1_H$, there exist an open set $V'\ni 1_G$ and an open neighborhood $U'$ of $x_0$ such that
$$\forall x\in(U'\cap C)\,\forall^*g\in V'\,(f(gx)\in W'f(x)).$$
\end{enumerate}

Let $W'\ni 1_H$ be any open subset of $H$. Then we find two open sets $V'\ni 1_G$ and $U'\ni x_0$ satisfying (iii).
Given any open set $O\subseteq Y$ with $f(x_0)\in O$ or $f(y_0)\in O$, without loss of generality, we can assume that $f(x_0)\in O$.
Since $f$ is continuous on $C$ and $x_0\in U'\cap C$, by shrinking, we can assume that $f(U'\cap C)\subseteq O$. Since $x_0\leftrightsquigarrow_V^\alpha y_0$, by Lemma~\ref{V^xV^y}, we have
$$\exists^*g^x\in V\,\exists^*g^y\in V\,(g^xx_0,g^yy_0\in U'\wedge\forall\beta<\alpha\,(g^xx_0\leftrightsquigarrow_{V'}^\beta g^yy_0)).$$
Moveover, since $x_0,y_0\in U\cap C$, by (1) and (ii), we can find $g^x,g^y\in V$ such that $g^xx_0,g^yy_0\in(U'\cap C)$, $g^xx_0\leftrightsquigarrow_{V'}^\beta g^yy_0~(\forall \beta<\alpha)$, $f(g^xx_0)\in Wf(x_0)$ and $f(g^yy_0)\in Wf(y_0)$. Thus there exist $h^x,h^y\in W$ so that $f(g^xx_0)=h^xf(x_0)$ and $f(g^yy_0)=h^yf(y_0)$.

Now we prove $f(x_0)\leftrightsquigarrow_W^\alpha f(y_0)$ by induction on $\alpha$. Since $g^xx_0,g^yy_0\in(U'\cap C)$, we have $h^xf(x_0),h^yf(y_0)\in O$. Thus we only need to show that $h^xf(x_0)\leftrightsquigarrow_{W'}^\beta h^yf(y_0)$ for $\beta<\alpha$.

For any open neighborhood $N$ of $h^xf(x_0)=f(g^xx_0)$, we have $g^xx_0\in f^{-1}(N)$. By the continuity of $f\upharpoonright C$, there is an open set $V_0\ni 1_G$ such that $(V_0g^xx_0\cap C)\subseteq f^{-1}(N)$. Since $g^yy_0\in(U'\cap C)$, we can find $v\in V_0\cap V'$ such that $vg^yy_0\in C$ and $f(vg^yy_0)\in W'f(g^yy_0)=W'h^yf(y_0)$. So $N\cap(W'h^yf(y_0))\ne\emptyset$. It follows that $h^xf(x_0)\in\overline{W'h^yf(y_0)}$. Similarly, we have $h^yf(y_0)\in\overline{W'h^xf(x_0)}$. Therefore, $h^xf(x_0)\leftrightsquigarrow_{W'}^0 h^yf(y_0)$.

If $\alpha=1$, we have already proved $f(x_0)\leftrightsquigarrow_W^1 f(y_0)$.

If $\alpha>1$, then for all $0<\beta<\alpha$, by inductive hypothesis, $g^xx_0\leftrightsquigarrow_{V'}^\beta g^yy_0$ implies $h^xf(x_0)=f(g^xx_0)\leftrightsquigarrow_{W'}^\beta f(g^yy_0)=h^yf(y_0)$. This completes the proof.
\hfill$\Box$ of Claim.

Put $V=G,W=H$, and $U=X$, then the proceeding Claim gives that, $x\leftrightsquigarrow_G^\alpha y$ implies $f(x)\leftrightsquigarrow_H^\alpha f(y)$ for $x,y\in C$ and $0<\alpha<\omega_1$.
\end{proof}

Now we return to $E(G)$ and the left multiplication action $c(G)\curvearrowright G^\omega$.

\begin{lemma}\label{right}
Let $G$ be a Polish group, $V\subseteq c(G)$ an open neighborhood of $1_{c(G)}$, $\alpha<\omega_1$, and $x,y,z\in G^\omega$. If $x\leftrightsquigarrow_V^\alpha y$, then $xz\leftrightsquigarrow_V^\alpha yz$.
\end{lemma}

\begin{proof}
We prove by induction on $\alpha$.

If $\alpha=0$, then we have $x\in\overline{Vy}$ and $y\in\overline{Vx}$. So $xz\in\overline{Vyz}$ and $yz\in\overline{Vxz}$, i.e., $xz\leftrightsquigarrow_V^0 yz$.

If $\alpha>0$, then for any open $W\ni 1_{c(G)}$ and any open $U\subseteq G^\omega$ with $xz\in U$ or $yz\in U$, we have $x\in Uz^{-1}$ or $y\in Uz^{-1}$, so there exist $g^x,g^y\in V$ such that $g^xx,g^yy\in Uz^{-1}$ and $g^xx\leftrightsquigarrow_W^\beta g^yy$ for $\beta<\alpha$. Then we have $g^xxz,g^yyz\in U$, and by induction hypothesis, $g^xxz\leftrightsquigarrow_W^\beta g^yyz$ for $\beta<\alpha$. Therefore, $xz\leftrightsquigarrow_V^\alpha yz$.
\end{proof}

\begin{definition}
Let $G$ be a Polish group, $0<\alpha<\omega_1$. We say $G$ is {\it $\alpha$-l.m.-unbalanced} if there exist $x,y\in G^\omega$ such that $x\leftrightsquigarrow_{c(G)}^\alpha y$ and $[x]\ne[y]$, or equivalently, there exists $z\in G^\omega\setminus c(G)$ such that $z\leftrightsquigarrow_{c(G)}^\alpha 1_{G^\omega}$.
\end{definition}

It is clear that $G$ is $0$-l.m.-unbalanced if $G$ contains at least two elements. And Lemma~\ref{1-unbalanced} gives that $G$ is l.m.-unbalanced iff it is $1$-l.m.-unbalanced.

\begin{theorem}\label{alpha-unbalanced}
Let $G,H$ be Polish groups, $0<\alpha<\omega_1$. If $G$ is $\alpha$-l.m.-unbalanced but $H$ is not, then $E(G)\not\le_B E(H)$.
\end{theorem}

\begin{proof}
Assume for contradiction that there is a Borel reduction $\theta$ of $E(G)$ to $E(H)$. From Theorem~\ref{homomorphism}, there is a comeager set $C\subseteq G^\omega$ such that for any $x,y\in C$, if $x\leftrightsquigarrow_{c(G)}^\alpha y$, then $\theta(x)\leftrightsquigarrow_{c(H)}^\alpha\theta(y)$. Since $G$ is $\alpha$-l.m.-unbalanced, there exists $z\in G^\omega\setminus c(G)$ such that $z\leftrightsquigarrow_{c(G)}^\alpha 1_{G^\omega}$. Since $C\cap zC\ne\emptyset$, we can find an $x\in(C\cap zC)$ and put $x=zy$. Then $x,y\in C$ and $xy^{-1}=z\not\in c(G)$, so $[x]\ne[y]$. By Lemma~\ref{right}, $x\leftrightsquigarrow_{c(G)}^\alpha y$. It follows that $\theta(x)\leftrightsquigarrow_{c(H)}^\alpha\theta(y)$ and $[\theta(x)]\ne[\theta(y)]$, contradicting that $H$ is not $\alpha$-l.m.-unbalanced.
\end{proof}

\subsection{Examples of $\alpha$-l.m.-unbalanced groups}

Let $G$ be a topological group, recall that $G$ is said to be \textit{distal} provided that $1_G\notin\overline{\{ghg^{-1}:g\in G\}}$ for all $h\ne 1_G\in G$ (cf.~\cite{rosenblatt}). It is trivial that all TSI Polish groups are distal. Moveover, if $H$ is a closed normal nilpotent subgroup of $G$ and $G/H$ is compact, then $G$ is distal (cf. \cite[\S 3. Proposition]{rosenblatt}). In particular, all nilpotent topological groups are distal.

\begin{theorem}\label{(*)}
Let $G$ be a Polish group. If $G$ is not distal, then $G$ is l.m.-unbalanced.
\end{theorem}

\begin{proof} Since $G$ is not distal, there exist $h\ne 1_G$ in $G$ and a sequence $(g_m)$ in $G$ with $g_mhg_m^{-1}\to 1_G$. For each $m,n\in\omega$, let $\gamma_m(n)=g_m$ and $\eta_m(n)=g_m^{-1}$. Then $\gamma_m,\eta_m\in c(G)$. We define $z\in G^\omega$ as
$$z(n)=\left\{\begin{array}{ll}1_G, & n=2k,\cr h, & n=2k+1.\end{array}\right.$$
It is clear that $\gamma_mz\eta_m\rightrightarrows 1_{G^\omega}$. From Lemma~\ref{unbalanced}, for any $x,y\in G^\omega$, if $xy^{-1}=z$, we have $x\leftrightsquigarrow y$. But it is trivial that $[x]\ne[y]$.
\end{proof}

\begin{theorem}\label{lc-nondistal}
Let $G$ be a locally compact Polish group. Then $G$ is not distal iff it is l.m.-unbalanced.
\end{theorem}

\begin{proof}
From Theorem~\ref{(*)}, we only need to prove that, if $G$ is l.m.-unbalanced then it is not distal.

Let $x,y\in G^\omega$ with $x\leftrightsquigarrow y$ but $[x]\ne[y]$. By Lemma~\ref{unbalanced}, there exist two sequences $(\gamma_m),(\eta_m)$ in $c(G)$ such that $\gamma_mxy^{-1}\eta_m\rightrightarrows 1_{G^\omega}$. Since $\gamma_m,\eta_m\in c(G)$ for each $m$, we can find $g_m,v_m\in G$ such that $g_m=\lim_n\gamma_m(n)$ and $v_m=\lim_n\eta_m(n)$.

Let $B\subseteq G$ be an open neighborhood of $1_G$ with $\overline{B}$ compact. There is a large enough $k\in\omega$ such that $\gamma_k(n)x(n)y(n)^{-1}\eta_k(n)\in B\subseteq\overline{B}$ for all $n\in\omega$. Thus there exist a $p\in\overline{B}$ and a strictly increasing natural numbers $n_0<n_1<\cdots$ such that $\lim_i\gamma_k(n_i)x(n_i)y(n_i)^{-1}\eta_k(n_i)=p$. By $[x]\ne[y]$, we have $xy^{-1}\notin c(G)$, thus $\gamma_kxy^{-1}\eta_k\notin c(G)$. So there exist a $q\in\overline{B}$ with $q\ne p$ and another strictly increasing natural numbers $l_0<l_1<\cdots$ such that $\lim_i\gamma_k(l_i)x(l_i)y(l_i)^{-1}\eta_k(l_i)=q$. Therefore, we have
$$\lim_ix(n_i)y(n_i)^{-1}=g_k^{-1}pv_k^{-1},\quad\lim_ix(l_i)y(l_i)^{-1}=g_k^{-1}qv_k^{-1}.$$

By the property of uniformly convergence of $\gamma_mxy^{-1}\eta_m\rightrightarrows 1_{G^\omega}$, we have
$$\begin{array}{ll}\lim_mg_mg_k^{-1}pv_k^{-1}v_m &=\lim_m\lim_i\gamma_m(n_i)x(n_i)y(n_i)^{-1}\eta_m(n_i)\cr
&=\lim_i\lim_m\gamma_m(n_i)x(n_i)y(n_i)^{-1}\eta_m(n_i)\cr
&=\lim_i1_G=1_G.\end{array}$$
Similarly, we have $\lim_mg_mg_k^{-1}qv_k^{-1}v_m=1_G$. In the end, we set $h=g_k^{-1}pq^{-1}g_k\ne 1_G$. It follows that
$$\lim_mg_mhg_m^{-1}=\lim_m(g_mg_k^{-1}pv_k^{-1}v_m)(g_mg_k^{-1}qv_k^{-1}v_m)^{-1}=1_G.$$
So $G$ is not distal.
\end{proof}

\begin{corollary}
Let $G,H$ be two Polish groups. If $G$ is not distal, and $H$ is TSI, or locally compact and distal, then $E(G)\not\le_BE(H)$.
\end{corollary}

\begin{proof}
It follows from theorems~\ref{TSI},~\ref{unbalanced to not},~\ref{(*)}, and~\ref{lc-nondistal}.
\end{proof}

\begin{example}
Let $G$ be the group of all $2\times 2$ real upper triangular matrices whose determinant $=1$, equipped with the usual topology. Then $h=\left(\begin{array}{cc}1 & 1\cr 0 & 1\end{array}\right)$ and $g_m=\left(\begin{array}{cc}1/m & 0\cr 0 & m\end{array}\right)$ witness that $G$ is not distal.
\end{example}

\begin{example}
Let $\Lambda$ be an infinite countable discrete group, $G$ a Polish group containing at least two elements, then the wreath product $\Lambda\wr G$ is not distal. To see this, fix an $a\in G$ with $a\ne 1_G$. We define $\chi:\Lambda\to G$ as $\chi(\lambda)=\left\{\begin{array}{ll}a, & \lambda=1_\Lambda,\cr 1_G, & \lambda\ne 1_\Lambda.\end{array}\right.$ and let $\Lambda=\{\lambda_m:m\in\omega\}$. Then $h=(1_\Lambda,\chi)$ and $g_m=(\lambda_m,1_{G^\Lambda})$ witness that $\Lambda\wr G$ is not distal.
\end{example}

\begin{theorem}
Let $G$ be a locally compact Polish group. Then $G$ is not $2$-l.m.-unbalanced.
\end{theorem}

\begin{proof}
Let $d$ be a compatible metric on $G$ and $d_u$ the supremum metric on $G^\omega$. Since $G$ is locally compact, we can find an $r>0$ with $\overline{B(1_G,r)}$ is compact. Let $V=\{\gamma\in c(G):d_u(1_{c(G)},\gamma)<r\}$.

Assume for contradiction that there exist $x,y\in G^\omega$ with $x\leftrightsquigarrow_{c(G)}^2 y$ but $[x]\ne[y]$, then there exist $g^x,g^y\in c(G)$ such that $g^xx\leftrightsquigarrow_V^1g^yy$. Following the arguments of $\Rightarrow$ part in the proof of Lemma~\ref{unbalanced}, we can find two sequences $(\gamma_m)$ and $(\eta_m)$ in $V$ such that $\gamma_mg^xxy^{-1}(g^y)^{-1}\eta_m\rightrightarrows 1_{G^\omega}$. Let $g_m=\lim_n\gamma_m(n)$. Following the arguments in the proof of Theorem~\ref{lc-nondistal}, there exists $h\ne 1_G$ such that $g_mhg_m^{-1}\to 1_G$.

Since $\gamma_m\in V$, we have $g_m\in\overline{B(1_G,r)}$. There exist a subsequence of $(g_m)$ converging to some $g\in G$, so $ghg^{-1}=1_G$, i.e., $h=1_G$. A contradiction!
\end{proof}

For $c(G)\curvearrowright G^\omega$, Lemma~\ref{unbalanced} can simplify arguments concerning $\leftrightsquigarrow$. Following the same spirit, we define a notion $\leftrightarrow_V^\alpha$ to simplify arguments concerning $\leftrightsquigarrow_V^\alpha$.

\begin{definition}\label{leftrightarrow}
Let $X$ be a Polish $G$-space, $V$ an open neighborhood of $1_G$, and let $\alpha<\omega_1$. We define $\leftrightarrow_V^\alpha$ by induction. We say
\begin{enumerate}
\item[(1)] $x\leftrightarrow_V^0y$, if $y\in\overline{Vx}$ and $x\in\overline{Vy}$;
\item[(2)] $x\leftrightarrow_V^\alpha y$ for $\alpha>0$, if for any open $W\ni 1_G$, there exist $g^x,g^y\in V$ such that $g^xx\leftrightarrow_W^\beta g^yy$ for all $\beta<\alpha$.
\end{enumerate}
\end{definition}

From their definitions, it is trivial to see that $x\leftrightsquigarrow_V^\alpha y$ implies $x\leftrightarrow_V^\alpha y$. Moreover, if $x\leftrightarrow_V^\alpha y$, then we have $x\leftrightarrow_V^{\alpha'} y$ for $0<\alpha'\le\alpha$, but $x\leftrightarrow_V^0 y$ may fails.

\begin{lemma}\label{equi-arrows}
Let $G$ be a Polish group, $c(G)\curvearrowright G^\omega$ the left multiplication action. Then $x\leftrightarrow_{c(G)}^\alpha y$ implies $x\leftrightsquigarrow_{c(G)}^\alpha y$ for all $x,y\in G^\omega$ and $\alpha<\omega_1$.
\end{lemma}

\begin{proof}
Let $d\leq 1$ be a compatible metric on $G$ and $d_u$ the supremum metric on $G^\omega$. For $\varepsilon>0$, we define $V_\varepsilon=\{\gamma\in c(G):d_u(1_{c(G)},\gamma)<\varepsilon\}$.

Let $x,y,x',y'\in G^\omega$, $m\in\omega$, and $\alpha<\omega_1$. If
\begin{enumerate}
\item[(1)] $x\leftrightarrow_{V_\varepsilon}^\alpha y$,
\item[(2)] $x(n)=x'(n),y(n)=y'(n)$ for $n\ge m$, and
\item[(3)] $x'(y')^{-1},y'(x')^{-1}\in\overline{V_{\varepsilon/2}}$,
\end{enumerate}
then we claim that $x'\leftrightsquigarrow_{V_\varepsilon}^\alpha y'$.

To see this, we prove by induction on $\alpha$. It is trivial that (3) implies $x'\leftrightsquigarrow_{V_\varepsilon}^0 y'$. For any $\delta>0$ and any open $U\subseteq G^\omega$ with $x'\in U$ or $y'\in U$, without loss of generality, we can assume that $x'\in U$, and $U=U_0\times\cdots\times U_l\times G^\omega$ with $l\ge m$ and $U_0,\cdots,U_l$ open in $G$. By $x\leftrightarrow_{V_\varepsilon}^\alpha y$, there exist $\gamma^x,\gamma^y\in V_{\varepsilon}$ such that $\gamma^xx\leftrightarrow_{V_{\delta/2}}^\beta\gamma^yy$ for $\beta<\alpha$. We put
$$(\gamma^x)'(n)=\left\{\begin{array}{ll}1_G, & n\le l,\cr \gamma^x(n), & n>l,\end{array}\right.
\quad(\gamma^y)'(n)=\left\{\begin{array}{ll}x'(n)y'(n)^{-1}, & n\le l,\cr \gamma^y(n), & n>l.\end{array}\right.$$
Then we have $(\gamma^x)',(\gamma^y)'\in V_\varepsilon$ and $(\gamma^x)'x',(\gamma^y)'y'\in U$. Note that $\gamma^xx\leftrightarrow_{V_{\delta/2}}^0\gamma^yy$ gives $\gamma^xx(\gamma^yy)^{-1},\gamma^yy(\gamma^xx)^{-1}\in\overline{V_{\delta/2}}$. Since $$(\gamma^x)'(n)x'(n)((\gamma^y)'(n)y'(n))^{-1}=\left\{\begin{array}{ll}1_G , & n\le l,\cr \gamma^x(n)x(n)(\gamma^y(n)y(n))^{-1}, & n>l,\end{array}\right.$$
we have $(\gamma^x)'x'((\gamma^y)'y')^{-1}\in\overline{V_{\delta/2}}$. Similarly, $(\gamma^y)'y'((\gamma^x)'x')^{-1}\in\overline{V_{\delta/2}}$. Since $\gamma^xx\leftrightarrow_{V_{\delta/2}}^\beta\gamma^yy$ implies $\gamma^xx\leftrightarrow_{V_\delta}^\beta\gamma^yy$ for $\beta<\alpha$. By the induction hypothesis, $(\gamma^x)'x'\leftrightsquigarrow_{V_\delta}^\beta(\gamma^y)'y'$ for $\beta<\alpha$. Thus $(\gamma^x)'$ and $(\gamma^y)'$ witnesses that $x'\leftrightsquigarrow_{V_\varepsilon}^\alpha y'$.

In the end, by setting $V_{\varepsilon/2}=c(G)$, it follows that $x\leftrightarrow_{c(G)}^\alpha y$ implies $x\leftrightsquigarrow_{c(G)}^\alpha y$.
\end{proof}

Given two sets $X,Y$ and a map $f:X\to Y$, recall that the map $f^\omega:X^\omega\to Y^\omega$ is defined as: $f^\omega(x)(n)=f(x(n))$ for $x\in X^\omega$ and $n\in\omega$.

\begin{lemma}\label{grouphomomorphism}
Let $G,H$ be two Polish groups, $\phi:G\to H$ a continuous homomorphism. Then for any $x,y\in G^\omega$, $\alpha<\omega_1$, open neighborhood $W\subseteq c(H)$ of $1_{c(H)}$,
and any open neighborhood $V\subseteq c(G)$ of $1_{c(G)}$ with $\phi^\omega(V)\subseteq W$, we have $x\leftrightarrow_V^\alpha y$ implies $\phi^\omega(x)\leftrightarrow_W^\alpha\phi^\omega(y)$.
\end{lemma}

\begin{proof}
Let $d$ be a compatible metric on $G$, and let $d_u$ be the supremum metric on $G^\omega$. By the continuity of $\phi$, we can see that $\phi^\omega:G^\omega\to H^\omega$ is continuous, and $\phi^\omega\upharpoonright c(G)$ is also continuous from $c(G)$ to $c(H)$.

We prove by induction on $\alpha$. If $\alpha=0$, note that $x\leftrightarrow_V^0 y$ iff $xy^{-1},yx^{-1}\in\overline V$. It follows from the continuity of $\phi^\omega$ that $\phi^\omega(x)\phi^\omega(y)^{-1}=\phi^\omega(xy^{-1})\in\overline W$ and $\phi^\omega(y)\phi^\omega(x)^{-1}=\phi^\omega(yx^{-1})\in\overline W$, and hence $\phi^\omega(x)\in\overline{W\phi^\omega(y)}$ and $\phi^\omega(y)\in\overline{W\phi^\omega(x)}$, i.e., $\phi^\omega(x)\leftrightarrow_W^0\phi^\omega(y)$.

For $\alpha>0$, we let $W'\subseteq c(H)$ be an open neighborhood of $1_{c(H)}$. By the continuity of $\phi^\omega\upharpoonright c(G)$, we can find an open neighborhood $V'$ of $1_{c(G)}$ such that $\phi^\omega(V')\subseteq W'$.
Note that $x\leftrightarrow_V^\alpha y$, so there exist $\gamma^x,\gamma^y\in V$ such that $\gamma^xx\leftrightarrow_{V'}^\beta\gamma^yy$ for $\beta<\alpha$. By induction hypothesis,
we have $\phi^\omega(\gamma^x)\phi^\omega(x)=\phi^\omega(\gamma^xx)\leftrightarrow_{W'}^\beta\phi^\omega(\gamma^yy)=\phi^\omega(\gamma^y)\phi^\omega(y)$.
Therefore, $\phi^\omega(\gamma^x)$ and $\phi^\omega(\gamma^y)$ witnesses that $\phi^\omega(x)\leftrightarrow_W^\alpha\phi^\omega(y)$.
\end{proof}

\begin{theorem}\label{alpha+1}
Let $G$ be a Polish group, $\Lambda$ an infinite countable discrete group, and $\alpha<\omega_1$. Then $G$ is $\alpha$-l.m.-unbalanced iff $\Lambda\wr G$ is $(\alpha+1)$-l.m.-unbalanced.
\end{theorem}

\begin{proof}
Let $\rho$ be a compatible metric on $\Lambda\wr G$, $\rho_u$ the supremum metric on $(\Lambda\wr G)^\omega$. For any $\varepsilon>0$, we define
$$V_\varepsilon=\{\gamma\in c(\Lambda\wr G):\rho_u(1_{c(\Lambda\wr G)},\gamma)<\varepsilon\}.$$

($\Rightarrow$). Suppose $G$ is $\alpha$-l.m.-unbalanced. Then there exists an $a\in G^\omega\setminus c(G)$ such that $a\leftrightsquigarrow_{c(G)}^\alpha 1_{G^\omega}$, and hence $a\leftrightarrow_{c(G)}^\alpha 1_{G^\omega}$. For any $n\in\omega$, let
$$\chi_n(\lambda)=\left\{\begin{array}{ll}a(n), & \lambda=1_\Lambda,\cr 1_G, & \lambda\ne 1_\Lambda.\end{array}\right.$$
We define $z\in(\Lambda\wr G)^\omega$ as $z(n)=(1_\Lambda,\chi_n)$ for $n\in\omega$. It is trivial that $z\notin c(\Lambda\wr G)$. To see that $\Lambda\wr G$ is $(\alpha+1)$-l.m.-unbalanced, by Lemma~\ref{equi-arrows}, we only need to show that $z\leftrightarrow_{c(\Lambda\wr G)}^{\alpha+1}1_{(\Lambda\wr G)^\omega}$.

Note that the underlying topology of $\Lambda\wr G$ is the product topology on $\Lambda\times G^\Lambda$. Thus for any given $\varepsilon>0$, there exists a $\lambda_0\in\Lambda$ such that, for any $\chi\in G^\Lambda$, if $\chi(\lambda)=1_G$ for all $\lambda\ne\lambda_0$, then $\rho(1_{\Lambda\wr G},(1_\Lambda,\chi))<\varepsilon/2$. Define
$$H=\{(1_\Lambda,\chi)\in\Lambda\wr G:\forall\lambda\in\Lambda\,(\lambda\ne\lambda_0\Rightarrow\chi(\lambda)=1_G)\}.$$
It is clear that $H$ is a closed subgroup of $\lambda\wr G$. Because $\rho_u(1_{(\Lambda\wr G)^\omega},x)\le\varepsilon/2<\varepsilon$ for all $x\in H^\omega$, we have
$$c(H)=(c(\Lambda\wr G)\cap H^\omega)\subseteq V_\varepsilon.$$
We define $f:G\to H$ as: for $g\in G$, $f(g)=(1_\Lambda,\chi^g)$ with
$$\chi^g(\lambda)=\left\{\begin{array}{ll}g, & \lambda=\lambda_0,\cr 1_G, & \lambda\ne\lambda_0.\end{array}\right.$$
It is easy to see that $f$ is a topological isomorphism from $G$ to $H$, and hence a continuous homomorphism from $G$ to $\Lambda\wr G$. Note that $f^\omega(c(G))\subseteq c(H)\subseteq V_\varepsilon$. Then Lemma~\ref{grouphomomorphism} gives $f^\omega(a)\leftrightarrow_{V_\varepsilon}^\alpha 1_{(\Lambda\wr G)^\omega}$.

In $\Lambda\wr G$, for $n\in\omega$, we have
$$(\lambda_0,1_{G^\Lambda})z(n)(\lambda_0,1_{G^\Lambda})^{-1}=(1_\Lambda,\chi^{a(n)})=f(a(n)).$$
Define $\gamma_0\in c(\Lambda\wr G)$ as $\gamma_0(n)=(\lambda_0,1_{G^\Lambda})$ for $n\in\omega$, then $\gamma_0z\gamma_0^{-1}=f^\omega(a)$.
It follows that $\gamma_0z=f^\omega(a)\gamma_0\leftrightarrow_{V_\varepsilon}^\alpha\gamma_0$. From Definition~\ref{leftrightarrow}, we can see that $\gamma_0z\leftrightarrow_{V_\varepsilon}^\beta\gamma_0$ for $0<\beta\le\alpha$. To prove for $\beta=0$, because $f^\omega(a),f^\omega(a)^{-1}\in H^\omega$, we have $\rho_u(1_{(\Lambda\wr G)^\omega},f^\omega(a))<\varepsilon$ and $\rho_u(1_{(\Lambda\wr G)^\omega},f^\omega(a)^{-1})<\varepsilon$. So $f^\omega(a),f^\omega(a)^{-1}\in\overline{V_\varepsilon}$, then we have $f^\omega(a)\gamma_0\in\overline{V_\varepsilon\gamma_0}$ and $\gamma_0\in\overline{V_\varepsilon f^\omega(a)\gamma_0}$, and hence $\gamma_0z=f^\omega(a)\gamma_0\leftrightarrow_{V_\varepsilon}^0\gamma_0$. It follows that $z\leftrightarrow_{c(\Lambda\wr G)}^{\alpha+1}1_{(\Lambda\wr G)^\omega}$.

($\Leftarrow$). On the other hand, suppose $\Lambda\wr G$ is $(\alpha+1)$-l.m.-unbalanced. Then there exists $z\in(\Lambda\wr G)^\omega\setminus c(\Lambda\wr G)$ such that $z\leftrightarrow_{c(\Lambda\wr G)}^{\alpha+1}1_{(\Lambda\wr G)^\omega}$. Thus for any $\varepsilon>0$, there exist $\gamma_\varepsilon^z,\gamma_\varepsilon\in c(\Lambda\wr G)$ such that $\gamma_\varepsilon^zz\leftrightarrow_{V_\varepsilon}^\alpha\gamma_\varepsilon$ and $\gamma_\varepsilon^zz\leftrightarrow_{V_\varepsilon}^0\gamma_\varepsilon$. Then we have $\gamma_\varepsilon^zz\gamma_\varepsilon^{-1}\leftrightarrow_{V_\varepsilon}^\alpha 1_{(\Lambda\wr G)^\omega}$ and $\gamma_\varepsilon^zz\gamma_\varepsilon^{-1}\in\overline{V_\varepsilon}$.

Fix a small enough $\varepsilon_0>0$ such that, for any $\gamma\in\overline{V_{\varepsilon_0}}$ and any $n\in\omega$, $\gamma(n)=(1_\Lambda,\chi)$ for some $\chi\in G^\Lambda$, i.e., $\overline{V_{\varepsilon_0}}$ is a subset of the subgroup $(\{1_\Lambda\}\times G^\Lambda)^\omega$ of $(\Lambda\wr G)^\omega$. So we can write $\gamma_{\varepsilon_0}^z(n)z(n)\gamma_{\varepsilon_0}(n)^{-1}=(1_\Lambda,\xi_n)$ with $\xi_n\in G^\Lambda$ for each $n\in\omega$. Note that $\gamma_{\varepsilon_0}^zz\gamma_{\varepsilon_0}^{-1}\notin c(\Lambda\wr G)$, we can find a $\lambda_0\in\Lambda$ such that $\lim_n\xi_n(\lambda_0)$ diverges.

Now we define $\pi:(\{1_\Lambda\}\times G^\Lambda)\to G$ as $\pi(1_\Lambda,\chi)=\chi(\lambda_0)$ for $\chi\in G^\Lambda$. It is clear that $\pi$ is a continuous homomorphism.
Note that $\gamma_{\varepsilon_0}^zz\gamma_{\varepsilon_0}^{-1}\in\overline{V_{\varepsilon_0}}\subseteq(\{1_\Lambda\}\times G^\Lambda)^\omega$. Then $\pi((\gamma_{\varepsilon_0}^zz\gamma_{\varepsilon_0}^{-1})(n))=\xi_n(\lambda_0)$ for all $n\in\omega$, so $\pi^\omega(\gamma_{\varepsilon_0}^zz\gamma_{\varepsilon_0}^{-1})\notin c(G)$. By the continuity of $\pi:(\{1_\Lambda\}\times G^\Lambda)\to G$, we have $\pi^\omega(c(\{1_\Lambda\}\times G^\Lambda))\subseteq c(G)$, and hence $\pi^\omega(V_{\varepsilon_0})\subseteq c(G)$. Then Lemma~\ref{grouphomomorphism} gives $\pi^\omega(\gamma_{\varepsilon_0}^zz\gamma_{\varepsilon_0}^{-1})\leftrightarrow_{c(G)}^\alpha 1_{G^\omega}$. Therefore, $G$ is $\alpha$-l.m.-unbalanced.
\end{proof}

\begin{corollary}
Let $G$ be a Polish group, $\Lambda$ an infinite countable discrete group, and $\alpha<\omega_1$. If $G$ is $\alpha$-l.m.-unbalanced but not $(\alpha+1)$-l.m.-unbalanced, then $E(G)<_BE(\Lambda\wr G)$.
\end{corollary}

\begin{proof}
It is trivial that $G$ is topologically isomorphic to a closed subgroup of $\Lambda\wr G$, so $E(G)\le_BE(\Lambda\wr G)$. It follows from Theorem~\ref{alpha-unbalanced} and Theorem~\ref{alpha+1} that $E(\Lambda\wr G)\not\le_BE(G)$.
\end{proof}

Let $G_0$ be a TSI Polish group containing at least two elements, and let $\Lambda$ be an infinite countable discrete group. We define $G_{n+1}=\Lambda\wr G_n$ for $n<\omega$. Then $G_n$ is $n$-l.m.-unbalanced but not $(n+1)$-l.m.-unbalanced. Therefore,
$$E(G_0)<_B\cdots<_BE(G_n)<_BE(G_{n+1})<_B\cdots.$$
From~\cite[Theorem 2.2.11]{gaobook}, all these $G_n$ are CLI. So far, we have not found any example of $\omega$-l.m.-unbalanced CLI Polish group.

On the other hand, if $G$ is a non-CLI Polish group, let $z\in G^\omega$ and $(\gamma_m)\in c(G)^\omega$ be which defined in the proof of Theorem~\ref{CLI}. We can inductively prove that, for any $V\ni 1_{c(G)}$, there exists $m$ large enough, such that $\gamma_m^{-1}z\leftrightarrow_V^\alpha 1_{G^\omega}$. So $z\leftrightarrow_{c(G)}^\alpha 1_{G^\omega}$, and hence $G$ is $\alpha$-l.m.-unbalanced for all ordinal $\alpha$ (even if $\alpha\ge\omega_1$).

\section{On TSI Polish groups}

\begin{definition}
Let $G$ be a Polish group. We define equivalence relation $E_*(G)$ on $G^\omega$ as: for $x,y\in G^\omega$,
$$xE_*(G)y\iff\lim_nx(0)x(1)\cdots x(n)y(n)^{-1}\cdots y(1)^{-1}y(0)^{-1}\mbox{ converges}.$$
\end{definition}

It is clear that $E(G)\sim_BE_*(G)$~(cf.~\cite[Proposition 2.2]{DZ}).

In this section, we focus on TSI Polish groups. It turns out that, for TSI Polish groups, $E_*(G)$ is a more convenient research object than $E(G)$.

\begin{lemma}
Let $G$ be a TSI Polish group, $d$ a complete compatible two-sided invariant metric on $G$.
\begin{enumerate}
\item[(1)] For $g_0,\cdots,g_n,h_0,\cdots,h_n\in G$, we have
$$d(g_0\cdots g_n,h_0\cdots h_n)=d(g_0\cdots g_nh_n^{-1}\cdots h_0^{-1},1_G)\le\sum_{k=0}^nd(g_k,h_k).$$
\item[(2)] For $x,y\in G^\omega$, we have
$$xE_*(G)y\iff\lim_n\sup_{n\le m}d(x(n)\cdots x(m),y(n)\cdots y(m))=0.$$
\item[(3)] For $x,y\in G^\omega$, if $xE_*(G)y$, then $\lim_nd(x(n),y(n))=0$.
\end{enumerate}
\end{lemma}

\begin{proof}
(1) Since $d$ is two-sided invariant, we have
$$\begin{array}{ll}d(g_0g_1,h_0h_1)&=d(g_0g_1h_1^{-1}h_0^{-1},1_G)=d(g_1h_1^{-1},g_0^{-1}h_0)\cr
& \le d(g_1h_1^{-1},1_G)+d(1_G,g_0^{-1}h_0)=d(g_0,h_0)+d(g_1,h_1).\end{array}$$
Then we can easily complete the proof of (1) by induction on $n$.

(2) By Cauchy criterion, $\lim_nx(0)x(1)\cdots x(n)y(n)^{-1}\cdots y(1)^{-1}y(0)^{-1}$ converges iff for any $\varepsilon>0$, there eixsts $N\in\omega$ such that, for $m>n>N$, we have
$$\begin{array}{ll}\varepsilon &>d(x(0)\cdots x(n)y(n)^{-1}\cdots y(0)^{-1},x(0)\cdots x(m)y(m)^{-1}\cdots y(0)^{-1})\cr
&=d(1_G,x(n+1)\cdots x(m)y(m)^{-1}\cdots y(n+1)^{-1})\cr
&=d(x(n+1)\cdots x(m),y(n+1)\cdots y(m)).\end{array}$$
It follows that
$$xE_*y\iff\lim_n\sup_{n\le m}d(x(n)\cdots x(m),y(n)\cdots y(m))=0.$$

(3) It is an easy corollary of (2).
\end{proof}

For the sake of brevity, we write
$$d(x,y)|_{[n,m+1)}=d(x,y)|_{[n,m]}=d(x(n)\cdots x(m),y(n)\cdots y(m)).$$

\subsection{Borel reducibility}

For any given metric space $(M,d)$, recall that $E(M;0)$ is an equivalence relation on $M^\omega$ (cf.~\cite[Definition 3.2]{ding12}) defined as
$$xE(M;0)y\iff\lim_nd(x(n),y(n))=0$$
for $x,y\in M^\omega$. If $G$ is a Polish group, then $E(G;0)$ is independent of the choice of left-invariant compatible metric $d$ on $G$, since $d(x(n),y(n))\to 0$ iff $x(n)^{-1}y(n)\to 1_G$.

\begin{lemma}\label{timesc_0}
Let $G$ be a TSI Polish group, then we have
$$E(G)\times E(G;0)\le_BE(G).$$
\end{lemma}

\begin{proof}
We define $\theta:G^\omega\times G^\omega\to G^\omega$ as
$$\theta(x,x')(n)=\left\{\begin{array}{ll}x(k), & n=2k,\cr x(k)x'(k), & n=2k+1.\end{array}\right.$$
Then $\theta$ witnesses that $E(G)\times E(G;0)\le_BE(G)$.
\end{proof}

\begin{lemma}\label{quotient}
Let $G,H,K$ be three TSI Polish groups. Suppose $\psi:G\to H$ and $\varphi:H\to K$ are continuous homomorphisms with $\varphi(\psi(G))=K$ such that $\ker(\varphi\circ\psi)$ is non-archimedean. Then $E(G)\le_BE(H)\times E(G;0)$.
\end{lemma}

\begin{proof}
By~\cite[Corollary 2.3.4]{gaobook}, $K\cong G/\ker(\varphi\circ\psi)$. Without loss of generality we may assume that $K=G/\ker(\varphi\circ\psi)$. Let $\phi=\varphi\circ\psi$. Then $\phi:G\to G/\ker(\varphi\circ\psi)$ is a continuous surjective homomorphism. Let $d_G,d_H$ be two-sided invariant complete compatible metrics on $G$ and $H$ respectively. Let
$$\begin{array}{ll}d_\phi(\phi(g),\phi(g'))&=\inf\{d_G(hg,h'g'):h,h'\in\ker(\varphi\circ\psi)\}\cr &=\inf\{d_G(hg,g'):h\in\ker(\varphi\circ\psi)\}.\end{array}$$
Then $d_\phi$ is a two-sided invariant complete compatible metric on $G/\ker(\varphi\circ\psi)$ (c.f.~\cite[Exercise 2.2.7]{gaobook}). It is clear that $d_\phi(\phi(g),\phi(g'))\le d_G(g,g')$ for $g,g'\in G$.

We only need to prove that $E_*(G)\le_BE_*(H)\times E(G;0)$. Define $\vartheta:G^\omega\to H^\omega\times G^\omega$ as, for $x\in G^\omega$ and $n\in\omega$,
$$\vartheta(x)(n)=(\psi(x(n)),x(n)).$$

For $x,y\in G^\omega$, if $xE_*(G)y$, then $\lim_k\sup_{k\le m}d_G(x,y)|_{[k,m]}=0$. So $\lim_nd_G(x(n),y(n))=0$. Since $d_G,d_H$ are two-sided invariant and $\psi$ is a continuous homomorphism, we have
$$\lim_k\sup_{k\le m}d_H(\psi(x(\cdot)),\psi(y(\cdot)))|_{[k,m]}=0.$$
Thus we have $(\vartheta(x),\vartheta(y))\in E_*(H)\times E(G;0)$.

On the other hand, if $(\vartheta(x),\vartheta(y))\in E_*(H)\times E(G;0)$, then we have $\lim_nd_G(x(n),y(n))=0$, and there exits $h_0\in H$ such that
$$\begin{array}{ll}& \lim_n\psi(x(0)\cdots x(n)y(n)^{-1}\cdots y(0)^{-1})\cr
=&\lim_n\psi(x(0))\cdots\psi(x(n))\psi(y(n))^{-1}\cdots\psi(y(0))^{-1}\cr
=&h_0.\end{array}$$
From $\varphi(h_0)\in K=\varphi(\psi(G))=\phi(G)$, there exists $g_0\in G$ with $\varphi(h_0)=\phi(g_0)$.
By continuity of $\varphi$ and $\phi=\varphi\circ\psi$, we have
$$\lim_n\phi(x(0)\dots x(n)y(n)^{-1}\dots y(0)^{-1})=\varphi(h_0)=\phi(g_0).$$

Since $\ker(\varphi\circ\psi)$ is non-archimedean, there exists a sequence of open subgroups $(M_l)$ of $\ker(\varphi\circ\psi)$ which forms a neighborhood base of $1_G$ in $\ker(\varphi\circ\psi)$. For any $\varepsilon>0$, there exist $l\in\omega$ and $0<\varepsilon'<\varepsilon$ such that
$$\{h\in\ker(\varphi\circ\psi):d_G(1_G,h)<3\varepsilon'\}\subseteq M_l\subseteq\{h\in\ker(\varphi\circ\psi):d_G(1_G,h)<\varepsilon\}.$$
Then there is an $N\in\omega$ such that
$$d_G(x(n),y(n))<\varepsilon',$$
$$d_\phi(\phi(x(0)\cdots x(n)y(n)^{-1}\cdots y(0)^{-1}),\phi(g_0))<\varepsilon'$$
for $n>N$. By the definition of $d_\phi$, there exist $h_n\in\ker(\varphi\circ\psi)$ for each $n>N$ such that
$$d_G(x(0)\cdots x(n)y(n)^{-1}\cdots y(0)^{-1},h_ng_0)<\varepsilon'.$$
Note that
$$\begin{array}{ll}& d_G(x(0)\cdots x(n)x(n+1)y(n+1)^{-1}y(n)^{-1}\cdots y(0)^{-1},\cr
&\quad\quad x(0)\cdots x(n)y(n)^{-1}\cdots y(0)^{-1})\cr
=& d_G(x(n+1),y(n+1))<\varepsilon'.\end{array}$$
Then $d_G(1_G,h_n^{-1}h_{n+1})=d_G(h_ng_0,h_{n+1}g_0)<3\varepsilon'$. So $h_n^{-1}h_{n+1}\in M_l$. It follows that $h_n^{-1}h_m$ is in the open subgroup $M_l$ for $m\ge n>N$. Thus $d_G(h_n,h_m)<\varepsilon$, and hence
$$\begin{array}{ll}& d_G(x(0)\cdots x(n)y(n)^{-1}\cdots y(0)^{-1},x(0)\cdots x(m)y(m)^{-1}\cdots y(0)^{-1})\cr
<& 2\varepsilon'+d_G(h_ng_0,h_mg_0)<3\varepsilon.\end{array}$$
This gives $xE_*(G)y$.

Therefore, $\vartheta$ witnesses that $E_*(G)\le_BE_*(H)\times E(G;0)$.
\end{proof}

\begin{theorem}\label{[0,1]}
Let $G,H,K$ be three TSI Polish groups. Suppose $\psi:G\to H$ and $\varphi:H\to K$ are continuous homomorphisms with $\varphi(\psi(G))=K$ such that $\ker(\varphi\circ\psi)$ is non-archimedean. If the interval $[0,1]$ embeds into $H$, then $E(G)\le_BE(H)$.
\end{theorem}

\begin{proof}
By~\cite[Theorem 3.4.(ii)]{ding12}, we have $E(G;0)\le_BE([0,1];0)$.

Let $f:[0,1]\to H$ be an embedding. By the uniformly continuity of $f$ and $f^{-1}:f([0,1])\to[0,1]$, it is trivial that $E([0,1];0)\sim_BE(f([0,1]);0)\le_BE(H;0)$.

Then lemmas~\ref{timesc_0} and~\ref{quotient} give $E(G)\le_BE(H)$.
\end{proof}

\begin{lemma}
Let $G,H$ be two TSI Polish groups, $G_c$ and $H_c$ open normal subgroups of $G$ and $H$ respectively, and let $T_G$ meets each coset of $G_c$ at exact one point. Let $\phi:G_c\to H_c$ be a topological group isomorphism, and $\theta:G\to H$ such that, for $u,v\in T_G$ and $g\in G_c$,
\begin{enumerate}
\item[(1)] $\phi(ugu^{-1})=\theta(u)\phi(g)\theta(u)^{-1}$,
\item[(2)] $\theta(ug)=\theta(u)\phi(g)$,
\item[(3)] if $u\ne v$, then $\theta(u)H_c\ne\theta(v)H_c$.
\end{enumerate}
Then $\theta^\omega$ is a continuous reduction of $E_*(G)$ to $E_*(H)$.
\end{lemma}

\begin{proof}
From (2), it is trivial to see that $\theta:G\to H$ is continuous, so $\theta^\omega$ is also continuous. we only need to check that $\theta^\omega$ is a reduction.

Let $x,y\in G^\omega$. We write, for brevity, $x(n)=u_ng_n$ and $y(n)=v_nh_n$ such that $u_n,v_n\in T_G$ and $g_n,h_n\in G_c$ for each $n$. Let $d_G$ and $d_H$ be two-sided invariant complete compatible metrics on $G$ and $H$ respectively.

If $xE_*(G)y$, then we have $\lim_k\sup_{k\le m}d_G(x,y)|_{[k,m]}=0$. In particular, we have $\lim_nd_G(x(n),y(n))=0$. Since $G_c$ is open, there exists an $N\in\omega$ such that,
for $n>N$, $x(n)$ and $y(n)$ are in the same coset of $G_c$, i.e., $u_n=v_n$. Now let $m\ge k>N$, we have
$$\begin{array}{cl}& x(k)x(k+1)\cdots x(m)\cr
=& u_kg_ku_{k+1}g_{k+1}\cdots u_mg_m\cr
=&(u_kg_ku_k^{-1})(u_ku_{k+1}g_{k+1}u_{k+1}^{-1}u_k^{-1})\cdots(u_k\cdots u_mg_mu_m^{-1}\cdots u_k^{-1})u_k\cdots u_m\cr
\stackrel{\rm Def}{=}& g_k'g_{k+1}'\cdots g_m'u_k\cdots u_m,\end{array}$$
$$\begin{array}{cl}& y(k)y(k+1)\cdots y(m)\cr
=& u_kh_ku_{k+1}h_{k+1}\cdots u_mh_m\cr
=&(u_kh_ku_k^{-1})(u_ku_{k+1}h_{k+1}u_{k+1}^{-1}u_k^{-1})\cdots(u_k\cdots u_mh_mu_m^{-1}\cdots u_k^{-1})u_k\cdots u_m\cr
\stackrel{\rm Def}{=}& h_k'h_{k+1}'\cdots h_m'u_k\cdots u_m.\end{array}$$
Note that
$$\phi(g_k')=\phi(u_kg_ku_k^{-1})=\theta(u_k)\phi(g_k)\theta(u_k)^{-1},$$
$$\begin{array}{ll}\phi(g_{k+1}')&=\phi(u_ku_{k+1}g_{k+1}u_{k+1}^{-1}u_k^{-1})\cr
&=\theta(u_k)\phi(u_{k+1}g_{k+1}u_{k+1}^{-1})\theta(u_k^{-1})\cr
&=\theta(u_k)\theta(u_{k+1})\phi(g_{k+1})\theta(u_{k+1})^{-1}\theta(u_k)^{-1},\end{array}$$
$$\begin{array}{ll}\phi(g_m')&=\phi(u_k\cdots u_mg_mu_m^{-1}\cdots u_k^{-1})\cr
&=\theta(u_k)\phi(u_{k+1}\cdots u_mg_mu_m^{-1}\cdots u_{k+1}^{-1})\theta(u_k^{-1})\cr
&=\cdots\cr
&=\theta(u_k)\cdots\theta(u_m)\phi(g_m)\theta(u_m)^{-1}\cdots\theta(u_k)^{-1}.\end{array}$$
Therefore, we have
$$\begin{array}{cl}& \theta(x(k))\theta(x(k+1))\cdots\theta(x(m))\cr
=&\theta(u_k)\phi(g_k)\theta(u_{k+1})\phi(g_{k+1})\cdots\theta(u_m)\phi(g_m)\cr
=&\phi(g_k')\phi(g_{k+1}')\cdots\phi(g_m')\theta(u_k)\cdots\theta(u_m).\end{array}$$
Similarly, we have
$$\begin{array}{cl}& \theta(y(k))\theta(y(k+1))\cdots\theta(y(m))\cr
=&\theta(u_k)\phi(h_k)\theta(u_{k+1})\phi(h_{k+1})\cdots\theta(u_m)\phi(h_m)\cr
=&\phi(h_k')\phi(h_{k+1}')\cdots\phi(h_m')\theta(u_k)\cdots\theta(u_m).\end{array}$$
It follows that from $\lim_k\sup_{k\le m}d_G(x,y)|_{[k,m]}=0$ that
$$d_G(g_k'\cdots g_m',h_k'\cdots h_m')=d_G(x(k)\cdots x(m),y(k)\cdots y(m))\to 0.$$
Since $\phi$ is a topological group isomorphism, we have
$$d_H(\theta^\omega(x),\theta^\omega(y))|_{[k,m]}=d_H(\phi(g_k')\cdots\phi(g_m'),\phi(h_k')\cdots\phi(h_m'))\to 0.$$
And hence $\theta^\omega(x)E_*(H)\theta^\omega(y)$.

On the other hand, if $\theta^\omega(x)E_*(H)\theta^\omega(y)$, we have $\lim_nd_H(\theta(x(n)),\theta(y(n)))=0$. By $\theta(x(n))=\theta(u_n)\phi(g_n)$ and $\theta(y(n))=\theta(v_n)\phi(h_n)$, we have $\theta(u_n)H_c=\theta(v_n)H_c$, i.e., $u_n=v_n$, for $n$ large enough. Then the similar arguments give $xE_*(G)y$.
\end{proof}

For $u\in G$, define a topological group automorphism $\iota_u:G_c\to G_c$ as $\iota_u(g)=ugu^{-1}$. We define
$${\rm Inn}_G(G_c)=\{\iota_u:u\in G\}.$$

\begin{theorem}\label{inn}
Let $G,H$ be two TSI Polish groups containing at least two elements, $G_c$ and $H_c$ open normal subgroups of $G$ and $H$ respectively. If there exists a topological group isomorphism $\phi:G_c\to H_c$ such that $\phi{\rm Inn}_G(G_c)\phi^{-1}\subseteq{\rm Inn}_H(H_c)$, then $E_*(G)\le_BE_*(H)$.
\end{theorem}

\begin{proof}
Corollary~\ref{timesE_0} gives $E_*(H\times\Z)\sim_BE_*(H)$, so we only need to show that $E_*(G)\le_BE_*(H\times\Z)$. Define a topological group isomorphism $\phi_0$ from $G_c$ to the open normal subgroup $H_c\times\{0\}$ of $H\times\Z$ as $\phi_0(g)=(\phi(g),0)$ for $g\in G_c$.

Fix a $T_G\subseteq G$ such that $T_G$ meats each coset of $G_c$ at exact one point. For each $u\in T_G$, since $\phi\iota_u\phi^{-1}\in{\rm Inn}_H(H_c)$, there exists some $w\in H$ such that $\phi\iota_u\phi^{-1}=\iota_w$. Thus $\phi_0\iota_u\phi_0^{-1}=\iota_{(w,n)}$ for all $n\in\Z$. Therefore, we can find $\theta:T_G\to(H\times\Z)$ such that, for $u,v\in T_G$,
\begin{enumerate}
\item[(i)] $\phi_0\iota_u\phi_0^{-1}=\iota_{\theta(u)}$,
\item[(ii)] if $u\ne v$, then $\theta(u)(H_c\times\{0\})\ne\theta(v)(H_c\times\{0\})$.
\end{enumerate}
Then we can extend $\theta$ to a map $G\to H$ satisfying clauses (1)--(3) in the proceeding lemma.
\end{proof}

\subsection{Borel irreducibility}

\begin{definition}[Farah~\cite{farah}]
\begin{enumerate}
\item[(1)] A map $\psi:\prod_nX_n\to\prod_nX_n'$ is \textit{additive} if there exist $0=l_0<l_1<\cdots<l_j<\cdots$ and maps $T_j:X_j\to\prod_{n\in[l_j,l_{j+1})}X_n'$ such that, for $x\in\prod_nX_n$,
    $$\psi(x)=T_0(x(0))^\smallfrown T_1(x(1))^\smallfrown T_2(x(2))^\smallfrown\cdots.$$
\item[(2)] Let $E$ and $F$ be equivalence relations on $\prod_nX_n$ and $\prod_nX_n'$ respectively, we say $E$ is \textit{additive reducible} to $F$, denoted by $E\le_A F$, if there exists an additive reduction of $E$ to $F$.
\end{enumerate}
\end{definition}

Let $E$ be an equivalence relation on $\prod_nX_n$, and let $I\subseteq\omega$ be infinite. Fix an element $w\in\prod_{n\notin I}X_n$. For $x\in\prod_{n\in I}X_n$, define $x\oplus w\in\prod_nX_n$ as: $(x\oplus w)(n)=x(n)$ for $n\in I$ and $(x\oplus w)(n)=w(n)$ for $n\notin I$. We define $E|_I^w$ on $\prod_{n\in I}X_n$ as: for $x,y\in\prod_{n\in I}X_n$,
$$xE|_I^wy\iff(x\oplus w)E(y\oplus w).$$

Let $(F_n)$ be a sequence of finite sets. A special equivalence relation $E_0(\prod_nF_n)$ defined as: for $x,y\in\prod_nF_n$,
$$xE_0(\prod_nF_n)y\iff\exists m\,\forall n>m\,(x(n)=y(n)).$$

The following lemma converts a Borel reduction to an additive reduction. This turns out to be a powerful tool to prove Borel irreducibility.

\begin{lemma}
Suppose $G$ is a TSI Polish group. Let $(F_n)$ be a sequence of finite sets, $E$ a Borel equivalence relation on $\prod_nF_n$ with $E_0(\prod_nF_n)\subseteq E$. If $E\le_BE_*(G)$, then there exist an infinite $I\subseteq\omega$ and a $w\in\prod_{n\notin I}F_n$ such that $E|_I^w\le_AE_*(G)$.
\end{lemma}

\begin{proof}
The following proof is a modification of the proof of \cite[Theorem 2.2]{DH}, claims (i)--(iii). We omit some similar arguments.

Assume that $\theta$ is a Borel reduction of $E$ to $E_*(G)$. Let $d$ be a two-sided invariant complete compatible metric on $G$. Following claims (i), (ii), and the arguments after Claim (ii) in the proof of~\cite[Theorem 2.2]{DH}, we construct two sequences of natural numbers $0=n_0<n_1<n_2<\cdots$ and $0=l_0<l_1<l_2<\cdots$, a sequence $(s_j)$ with $s_j\in\prod_{n_j<n<n_{j+1}}F_n$ for each $j\in\omega$. Put $I=\{n_j:j\in\omega\}$, and put $w=\bigcup_js_j\in\prod_{n\notin I}F_n$. Our construction confirms that, for any $x,y\in\prod_{n\in I}F_n$, we have:
\begin{enumerate}
\item[(a)] if $x(n)=y(n)$ for $n>n_j$, then for $l_{j+1}\le k\le m$,
$$d(\theta(x\oplus w),\theta(y\oplus w))|_{[k,m]}<2^{-j};$$
\item[(b)] if $x(n)=y(n)$ for $n\le n_j$, then for $k\le m<l_{j+1}$,
$$d(\theta(x\oplus w),\theta(y\oplus w))|_{[k,m]}<2^{-j}.$$
\end{enumerate}

For each $n\in I$, fix an $a_n^\#\in F_n$. Define $p_j:F_{n_j}\to\prod_{n\in I}F_n$ for each $j\in\omega$ as $p_j(a)(n)=\left\{\begin{array}{ll}a, & n=n_j,\cr a_n^\#, & n\ne n_j\end{array}\right.$ for $a\in F_{n_j}$ and $n\in\omega$. Then we define $T_{n_j}:F_{n_j}\to G^{l_{j+1}-l_j}$ as: for $a\in F_{n_j}$ and $i<l_{j+1}-l_j$,
$$T_{n_j}(a)(i)=\theta(p_j(a)\oplus w)(l_j+i).$$
The additive mapping $\psi:\prod_{n\in I}F_n\to G^\omega$ is defined as: for $x\in\prod_{n\in I}F_n$,
$$\psi(x)=T_{n_0}(x(n_0))^\smallfrown T_{n_1}(x(n_1))^\smallfrown T_{n_2}(x(n_2))^\smallfrown\cdots.$$

Now we will show that $\psi$ is a reduction of $E|_I^w$ to $E_*(G)$.

For any $x\in\prod_{n\in I}F_n$ and $j\in\omega$, define $e_j(x),e'_j(x)\in\prod_{n\in I}F_n$ as
$$e_j(x)(n)=\left\{\begin{array}{ll}x(n), & n=n_j,\cr a_n^\#, & n\ne n_j,\end{array}\right.\quad
e'_j(x)(n)=\left\{\begin{array}{ll}x(n), & n\le n_j,\cr a_n^\#, & n>n_j.\end{array}\right.$$
Applying (a) for $l_j\le k\le m$, we have
$$d(\theta(e_j(x)\oplus w),\theta(e'_j(x)\oplus w))|_{[k,m]}<2^{-(j-1)}.$$
and applying (b) for $k\le m<l_{j+1}$, we have
$$d(\theta(x\oplus w)),\theta(e'_j(x)\oplus w))|_{[k,m]}<2^{-j};$$

{\bf Claim.} $\theta(x\oplus w)E_*(G)\psi(x)$.

\textit{Proof of Claim.}
Since $d$ is complete and two-sided invariant, by Cauchy criterion, we only need to show that
$$\lim_k\sup_{k\le m}d(\theta(x\oplus w),\psi(x))|_{[k,m]}\to 0.$$
For any $k\le m$, there exists $j\le j'$ with $l_j\le k<l_{j+1}$ and $l_{j'}\le m<l_{j'+1}$. If $j=j'$, note that $p_j(x(n_j))=e_j(x)$, so we have
$$\psi(x)\upharpoonright[l_j,l_{j+1})=\theta(e_j(x)\oplus w)\upharpoonright[l_j,l_{j+1}),$$
and hence
$$\begin{array}{ll}& d(\theta(x\oplus w),\psi(x))|_{[k,m]}\cr
\le & d(\theta(x\oplus w),\theta(e'_j(x)\oplus w))|_{[k,m]}+d(\theta(e'_j(x)\oplus w),\theta(e_j(x)\oplus w))|_{[k,m]}\cr
<& 2^{-j}+2^{-(j-1)}=3\cdot 2^{-j};\end{array}$$
otherwise, $j<j'$, we have
$$\begin{array}{ll}& d(\theta(x\oplus w),\psi(x))|_{[k,m]}\cr
\le &d(\theta(x\oplus w),\psi(x))|_{[k,l_{j+1})}+\sum_{i=j+1}^{j'-1}d(\theta(x\oplus w),\psi(x))|_{[l_i,l_{i+1})}\cr
&+d(\theta(x\oplus w),\psi(x))|_{[l_{j'},m]}\cr
< & 3\cdot\sum_{i=j}^{j'}2^{-i}<3\cdot 2^{-(j-1)}.\end{array}$$
This gives $\theta(x\oplus w)E_*(G)\psi(x)$.
\hfill$\Box$ of Claim.

In the end, for $x,y\in\prod_{n\in I}F_n$, $xE|_I^wy$ means $(x\oplus w)E(y\oplus w)$, so
$$xE|_I^wy\iff\theta(x\oplus w)E_*(G)\theta(y\oplus w)\iff\psi(x)E_*(G)\psi(y).$$
This competes the proof.
\end{proof}

Let $G,H$ be two TSI Polish groups, and let $d_G$ and $d_H$ be complete compatible two-sided invariant metrics on $G$ and $H$ respectively. Assume that $\theta$ is a Borel reduction of  $E_*(G)$ to $E_*(H)$. Let $(F_n)$ be a sequence of finite subsets of $G$ such that
\begin{enumerate}
\item[(i)] $1_G\in F_n=F_n^{-1}$,
\item[(ii)] $F_{n-1}^{3n+2}\subseteq F_n$, and
\item[(iii)] $\bigcup_nF_n$ is dense in $G$.
\end{enumerate}
Let $E$ be the restriction of $E_*(G)$ on $\prod_nF_n$. There exist an infinite $I\subseteq\omega$ and a $w\in\prod_{n\notin I}F_n$ such that $E|_I^w\le_AE_*(H)$.
So there are natural numbers $0=n_0<n_1<n_2<\cdots$ with $I=\{n_j:j\in\omega\}$, $0=l_0<l_1<l_2<\cdots$, $T_{n_j}:F_{n_j}\to H^{l_{j+1}-l_j}$, and $\psi:\prod_{n\in I}F_n\to H^\omega$ with
$$\psi(x)=T_{n_0}(x(n_0))^\smallfrown T_{n_1}(x(n_1))^\smallfrown T_{n_2}(x(n_2))^\smallfrown\cdots,$$
such that $\psi$ is an additive reduction of $E|_I^w$ to $E_*(H)$.

For $s=(h_0,\cdots,h_{l-1})$ and $t=(h'_0,\cdots,h'_{l-1})$ in $H^l$, let
$$d_H^\infty(s,t)=\max_{0\le k\le m<l}d_H(h_k\cdots h_m,h'_k\cdots h'_m).$$

\begin{lemma}
For any $q\in\omega$, there exists a $\delta_q>0$ such that
$$\forall^\infty n\in I\,\forall g,g'\in F_n\,(d_G(g,g')<\delta_q\Rightarrow d_H^\infty(T_n(g),T_n(g'))<2^{-q}).$$
\end{lemma}

\begin{proof}
If not, then there exist an $\varepsilon_0>0$, a strictly increasing sequence $(j_p)$, and $g_p,g'_p\in F_{n_{j_p}}$ for each $p$, such that
$$d_G(g_p,g'_p)<2^{-p},\quad d_H^\infty(T_{n_{j_p}}(g_p),T_{n_{j_p}}(g'_p))\ge\varepsilon_0.$$
For each $n\in I$, put
$$x(n)=\left\{\begin{array}{ll}g_p, & n=n_{j_p},\cr 1_G, & \mbox{otherwise},\end{array}\right.\quad
y(n)=\left\{\begin{array}{ll}g'_p, & n=n_{j_p},\cr 1_G, & \mbox{otherwise}.\end{array}\right.$$
Then for any $k\le m$, letting $n_{j_p}\le k<n_{j_{p+1}}$, we have
$$\begin{array}{ll}& d_G(x\oplus w,y\oplus w)|_{[k,m]}\cr
\le &\sum_{k\le n\le m}d_G((x\oplus w)(n),(y\oplus w)(n))\cr
\le &\sum_{r\ge p}d_G(g_r,g'_r)<\sum_{r\ge p}2^{-r}=2^{-(p-1)}.\end{array}$$
It follows that $(x\oplus w)E_*(G)(y\oplus w)$. But for any $p\in\omega$,
$$\max_{l_{j_p}\le k\le m<l_{j_{p+1}}}d_H(\psi(x),\psi(y))|_{[k,m]}
=d_H^\infty(T_{n_{j_p}}(g_p),T_{n_{j_p}}(g'_p))\ge\varepsilon_0.$$
So $(\psi(x),\psi(y))\notin E_*(H)$, contradicting that $\psi$ is a reduction.
\end{proof}

In the rest of this subsection, we assume that $H$ is locally compact.

Then there exists a $q_c\in\omega$ such that the closure of $W_c=\{h\in H:d_H(1_H,h)<2^{-q_c}\}$ is compact. Note that $W_c=W_c^{-1}$ and $hW_ch^{-1}=W_c$ for all $h\in H$.
Define $H_c=\bigcup_mW_c^m$. Then $H_c$ is an open normal subgroup of $H$. Let $V_c=\{g\in G:d_G(1_G,g)<\delta_{q_c}\}$ and $G_c=\bigcup_mV_c^m$.
Note that $V_c=V_c^{-1}$ and $gV_cg^{-1}=V_c$ for all $g\in G$. So $G_c$ is an open normal subgroup of $G$.

For $j>0$, we define
$$u_{n_j}=w(n_{j-1}+1)\cdots w(n_j-1),$$
and $u_{n_0}=u_0=1_G$.
For any $g\in\bigcup_nF_n$, if $g\in F_n$ for some $n<n_j$, we have $u_{n_j}^{-1}g\in F_{n_j-1}^{n_j}\subseteq F_{n_j}$. We define
$$S_{n_j}(g)=T_{n_j}(u_{n_j}^{-1}g)(0)\cdots T_{n_j}(u_{n_j}^{-1}g)(l_{j+1}-l_j-1).$$
Note that $d_G(u_{n_j}^{-1}g,u_{n_j}^{-1}g')=d_G(g,g')$. From the preceding lemma, we have
$$\forall q\,\forall^\infty j\,\forall g,g'\in F_{n_j}\,(d_G(g,g')<\delta_q\Rightarrow d_H(S_{n_j}(g),S_{n_j}(g'))<2^{-q}).$$

Fix a coset $C$ of $G_c$, note that $C$ is open, we can fix a $v\in C\cap\bigcup_nF_n$. For $g\in G_c\cap\bigcup_nF_n$, there exists an $m>0$ such that $g\in V_c^m$.
So we can find $1_G=g_0',g_1',\cdots,g_m'=g\in G_c$ such that $(g_i')^{-1}g_{i+1}'\in V_c$. Put $g_0=g_0'$ and $g_m=g_m'$. For $1\le i<m$, we can also find $g_i$ in some small enough neighborhood of $g_i'$ such that  $g_i\in G_c\cap\bigcup_nF_n$ and $g_i^{-1}g_{i+1}\in V_c$. Then $d_G(g_i,g_{i+1})=d_G(1_G,g_ig_{i+1}^{-1})<\delta_{q_c}$ for $i<m$. Moveover, $d_G(v^{-j}g_iv^{j+1},v^{-j}g_{i+1}v^{j+1})=d_G(g_i,g_{i+1})<\delta_{q_c}$ for all $j\in\omega$. For $j$ large enough, we have $v,g_i\in F_{n_j-1}$, and hence $v^{-j}g_iv^{j+1}\in F_{n_j-1}^{2j+2}\subseteq F_{n_j}$. So, for $i<m$,
$$d_H(S_{n_j}(v^{-j}g_i v^{j+1}),S_{n_j}(v^{-j}g_{i+1}v^{j+1}))<2^{-q_c}.$$
Thus $S_{n_j}(v^{-j}g v^{j+1})S_{n_j}(v)^{-1}\in W_c^m\subseteq\overline{W_c^m}$. We define
$$S^C_{n_j}(g)=S_{n_0}(v)\cdots S_{n_{j-1}}(v)S_{n_j}(v^{-j}gv^{j+1})S_{n_j}(v)^{-1}S_{n_{j-1}}(v)^{-1}\cdots S_{n_0}(v)^{-1}.$$
Then $S^C_{n_j}(g)\in W_c^m\subseteq\overline{W_c^m}$ and $S^C_{n_j}(1_G)=1_H$. Fix an ultrafilter $\mathfrak A$ on $\omega$ such that $\mathfrak A$ does not contain any finite set. Note that $\overline{W_c^m}=\overline W_c^m$ is compact. We now define
$$S^C(g)=\lim_{j\in\mathfrak A}S^C_{n_j}(g)\in\overline{W_c^m}.$$
For the definition and the existence of the ultrafilter limit $\lim_{j\in\mathfrak A}$, one can see~\cite[\S 1.5]{CK}. Since $\overline{W_c^m}\subseteq W_c^{m+1}\subseteq H_c$ for each $m\in\omega$, we have $S^C(g)\in H_c$. Note that
$$\forall q\,\forall g,g'\in G_c\cap\bigcup_nF_n\,(d_G(g,g')<\delta_q\Rightarrow d_H(S^C(g),S^C(g'))\le 2^{-q}).$$
So $S^C$ is uniformly continuous on $G_c\cap\bigcup_nF_n$. Therefore, we can extend $S^C$ to a uniformly continuous map from $G_c$ to $H_c$, which is still denoted as $S^C$ for brevity.

It is worth noting that $S^C(1_G)=1_H$.

For $x\in G_c^\omega$, recall that $(S^C)^\omega(x)\in H_c^\omega$ is defined as:
$$(S^C)^\omega(x)(p)=S^C(x(p))\quad(\forall p\in\omega).$$

\begin{lemma}\label{S^C}
For $x,y\in G_c^\omega$, if $\lim_pd_G(x(p),y(p))=0$, then
$$xE_*(G)y\iff (S^C)^\omega(x)E_*(H)(S^C)^\omega(y).$$
\end{lemma}

\begin{proof}
For each $p\in\omega$, we can find $g_p,g'_p\in G_c\cap\bigcup_nF_n$ such that
$$d_G(g_p,x(p))<2^{-p},\quad d_G(g'_p,y(p))<2^{-p},$$
$$d_H(S^C(g_p),S^C(x(p)))<2^{-p},\quad d_H(S^C(g'_p),S^C(y(p)))<2^{-p}.$$
We can find a large enough $j(p)\in\omega$ for $p\in\omega$ such that
\begin{enumerate}
\item[(0)] $j(0)<j(1)<j(2)<\cdots$;
\item[(1)] $v\in F_n$ for some $n<n_{j(0)}$;
\item[(2)] $g_p,g'_p\in F_n$ for some $n<n_{j(p)}$;
\item[(3)] for brevity, we write $n(p)=n_{j(p)}$, then
$$d_H(S^C(g_p),S^C_{n(p)}(g_p))<2^{-p},\quad d_H(S^C(g'_p),S^C_{n(p)}(g'_p))<2^{-p}.$$
\end{enumerate}

For $n\in I$, we define
$$\hat x(n)=\left\{\begin{array}{ll}u_{n(p)}^{-1}v^{-j(p)}g_pv^{j(p)+1}, & n=n(p),\cr u_n^{-1}v, & \mbox{otherwise,}\end{array}\right.$$
$$\hat y(n)=\left\{\begin{array}{ll}u_{n(p)}^{-1}v^{-j(p)}g'_pv^{j(p)+1}, & n=n(p),\cr u_n^{-1}v, & \mbox{otherwise.}\end{array}\right.$$
Note that $\hat x(n),\hat y(n)\in F_{n-1}^{3n+2}\subseteq F_n$ for each $n\in I$.

Then for $k\le m$, if $n(q-1)<k\leq m\leq n(q)$, we have
$$\begin{array}{ll} & d_G(\hat x\oplus w,\hat y\oplus w)|_{[k,m]}\cr
\leq&d_G(u_{n(q)}^{-1}v^{-j(q)}g_qv^{j(q)+1},u_{n(q)}^{-1}v^{-j(q)}g'_qv^{j(q)+1})\cr
=&d_G(g_q,g_q')\to 0\quad(q\to\infty).
\end{array}$$
If $n(q-1)<k\le n(q)$ and $n(r)\le m<n(r+1)~(q\leq r)$, we have
$$\begin{array}{ll} & d_G(\hat x\oplus w,\hat y\oplus w)|_{[k,m]}\cr
=& d_G(u_{n(q)}^{-1}v^{-j(q)}g_q\cdots g_{r}v^{j(r)+1},u_{n(q)}^{-1}v^{-j(q)}g'_q\cdots g'_{r}v^{j(r)+1})\cr
=& d_G(g_q\cdots g_{r},g'_q\cdots g'_{r}).\end{array}$$
Thus we have
$$\begin{array}{ll} & \left|d_G(\hat x\oplus w,\hat y\oplus w)|_{[k,m]}-d_G(x,y)|_{[q,r]}\right|\cr
=& \left|d_G(g_q\cdots g_r,g'_q\cdots g'_r)-d_G(x(q)\cdots x(r),y(q)\cdots y(r))\right|\cr
\le & d_G(g_q\cdots g_r,x(q)\cdots x(r))+d_G(g'_q\cdots g'_r,y(q)\cdots y(r))\cr
\le & \sum_{q\le p\le r}(d_G(g_p,x(p))+d_G(g'_p,y(p)))\cr
< & 2\cdot\sum_{q\le p}2^{-p}=2^{-(q-2)}\to 0\quad(q\to\infty).\end{array}$$
This implies that
$$xE_*(G)y\iff(\hat x\oplus w)E_*(G)(\hat y\oplus w).$$
Since $\psi$ is an additive reduction, we have
$$(\hat x\oplus w)E_*(G)(\hat y\oplus w)\iff\hat xE|_I^w\hat y\iff\psi(\hat x)E_*(H)\psi(\hat y).$$

On the other hand, for $k\le m$, if $l_j\le k<l_{j+1}$ and $l_{s}\le m<l_{s+1}$, then $j\le s$, so there exist $q$ and $r$ such that $n(q)\le n_j<n(q+1)$ and $n(r)<n_{s}\le n(r+1)$. Then
$$\begin{array}{ll} & d_H(\psi(\hat x),\psi(\hat y))|_{[k,m]}\cr
\le & d_H(\psi(\hat x),\psi(\hat y))|_{[k,l_{j+1})}+d_H(\psi(\hat x),\psi(\hat y))|_{[l_{j+1},l_{s})}+d_H(\psi(\hat x),\psi(\hat y))|_{[l_{s},m]}.\end{array}$$

Note that
$$d_H(\psi(\hat x),\psi(\hat y))|_{[k,l_{j+1})}\le d_H^\infty(T_{n_j}(\hat x(n_j)),T_{n_j}(\hat y(n_j))).$$
If $n_j\ne n(q)$, then $\hat x(n_j)=\hat y(n_j)$; if $n_j=n(q)$, since
$$\begin{array}{ll}d_G(\hat x(n(q)),\hat y(n(q))) & =d_G(g_q,g'_q)\cr
&\le d_G(x(q),y(q))+2^{-q}+2^{-q}\to 0\quad(q\to\infty),\end{array}$$
we have
$$d_H^\infty(T_{n(q)}(\hat x(n(q))),T_{n(q)}(\hat y(n(q))))\to 0\quad(q\to\infty).$$
Therefore, $d_H(\psi(\hat x),\psi(\hat y))|_{[k,l_{j+1})}\to 0$ as $k\to\infty$.

Similarly, we also have $d_H(\psi(\hat x),\psi(\hat y))|_{[l_{s},m]}\to 0$ as $m\to\infty$.

In the end, we have
$$\psi(\hat x)(l_{j+1})\cdots\psi(\hat x)(l_{s}-1)=S_{n_{j+1}}(u_{n_{j+1}}\hat x(n_{j+1}))\cdots S_{n_{s-1}}(u_{n_{s-1}}\hat x(n_{s-1})).$$
For $n_i\in I$ with $j+1\le i\le s-1$, if $n_i=n(p)$ for some $q+1\le p\le r$, then
$$S_{n_i}(u_{n_i}\hat x(n_i))=S_{n_{i-1}}(v)^{-1}\cdots S_{n_0}(v)^{-1}S^C_{n(p)}(g_p)S_{n_0}(v)\cdots S_{n_i}(v);$$
and if $n_i\ne n(p)$ for all $q+1\le p\le r$, then
$$S_{n_i}(u_{n_i}\hat x(n_i))=S_{n_i}(v).$$
This gives that
$$\begin{array}{ll}& \psi(\hat x)(l_{j+1})\cdots\psi(\hat x)(l_{s}-1)\cr
=& S_{n_j}(v)^{-1}\cdots S_{n_0}(v)^{-1}S^C_{n(q+1)}(g_{q+1})\cdots S^C_{n(r)}(g_r)S_{n_0}(v)\cdots S_{n_{s-1}}(v).\end{array}$$
Similarly, we have
$$\begin{array}{ll}& \psi(\hat y)(l_{j+1})\cdots\psi(\hat y)(l_{s}-1)\cr
=& S_{n_j}(v)^{-1}\cdots S_{n_0}(v)^{-1}S^C_{n(q+1)}(g'_{q+1})\cdots S^C_{n(r)}(g'_r)S_{n_0}(v)\cdots S_{n_{s-1}}(v).\end{array}$$
Thus we have
$$\begin{array}{ll} &d_H(\psi(\hat x),\psi(\hat y))|_{[l_{j+1},l_{s})}\cr
=& d_H(S^C_{n(q+1)}(g_{q+1})\cdots S^C_{n(r)}(g_r),S^C_{n(q+1)}(g'_{q+1})\cdots S^C_{n(r)}(g'_r)).\end{array}$$
Note that
$$\begin{array}{ll} & d_H(S^C_{n(q+1)}(g_{q+1})\cdots S^C_{n(r)}(g_r),S^C(x(q+1))\cdots S^C(x(r)))\cr
\le & \sum_{q+1\le p\le r}d_H(S^C_{n(p)}(g_p),S^C(x(p)))\cr
\le & \sum_{q+1\le p\le r}(d_H(S^C_{n(p)}(g_p),S^C(g_p))+d_H(S^C(g_p),S^C(x(p))))\cr
< & 2\cdot\sum_{q+1\le p\le r}2^{-p}<2^{-(q-1)}\to 0\quad(q\to\infty).\end{array}$$
Similarly, we also have
$$d_H(S^C_{n(q+1)}(g'_{q+1})\cdots S^C_{n(r)}(g'_r),S^C(y(q+1))\cdots S^C(y(r)))\to 0\quad(q\to\infty).$$

To sum up, we have
$$\begin{array}{ll}&\lim_kd_H(\psi(\hat x),\psi(\hat y))|_{[k,m]}=0\cr
\iff &\lim_jd_H(\psi(\hat x),\psi(\hat y))|_{[l_{j+1},l_{s})}=0\cr
\iff &\lim_qd_H((S^C)^\omega(x),(S^C)^\omega(y))|_{[q+1,r]}=0.\end{array}$$
It follows that
$$\psi(\hat x)E_*(H)\psi(\hat y)\iff (S^C)^\omega(x)E_*(H)(S^C)^\omega(y),$$
and hence $xE_*(G)y\iff (S^C)^\omega(x)E_*(H)(S^C)^\omega(y)$.
\end{proof}

Recall that a topological group $G$ is said to {\it have no small subgroups} if there exists an open set $U\ni 1_G$ such that $\{1_G\}$ is the only subgroup of $G$ contained in $U$ (cf. \cite[Definition 2.16]{HM07}).

\begin{lemma}\label{locally injective}
Suppose $G$ has no small subgroups. Then there is a $D>0$ such that, for all $1_G\ne h\in G_c$, if $d_G(1_G,h)\le D$, then
$$\inf_{g\in G_c}d_H(S^C(g),S^C(hg))>0.$$

In particular, for $g,g'\in G_c$, if $0<d_G(g,g')\le D$, then $S^C(g)\ne S^C(g')$.
\end{lemma}

\begin{proof}
If not, there exists a sequence $(h_q)$ in $G_c$ such that $h_q\ne 1_G$ and $\inf_{g\in G_c}d_H(S^C(g),S^C(h_qg))=0$ for each $q\in\omega$, and $\lim_qd_G(1_G,h_q)=0$.

Since $G$ has no small subgroups, there exists some $D_0>0$ such that, for each $q\in\omega$, we can find an $m_q>0$ such that $d_G(1_G,h_q^{m_q})\ge D_0$.  And for $i<m_q$, we can find an element $g_{q,i}\in G_c$ such that
$$d_H(S^C(g_{q,i}),S^C(h_qg_{q,i}))<2^{-(q+i+3)}.$$
Let $M_{-1}=0$ and $M_q=m_0+\cdots m_q$ for $q\in\omega$.

Now for $p\in\omega$, define
$$x(p)=\left\{\begin{array}{ll}g_{q,i}, & p=2(M_{q-1}+i),i<m_q,\cr g_{q,i}^{-1}, & p=2(M_{q-1}+i)+1,i<m_q,\end{array}\right.$$
$$y(p)=\left\{\begin{array}{ll}h_qg_{q,i}, & p=2(M_{q-1}+i),i<m_q,\cr g_{q,i}^{-1}, & p=2(M_{q-1}+i)+1,i<m_q.\end{array}\right.$$
By $\lim_qd_G(1_G,h_q)=0$, we have $\lim_pd_G(x(p),y(p))=0$, then it follows from Lemma~\ref{S^C} that
$$xE_*(G)y\iff (S^C)^\omega(x)E_*(H)(S^C)^\omega(y).$$

Note that $d_G(x,y)|_{[2M_{q-1},2M_q)}=d_G(1_G,h_q^{m_q})\ge D_0$, so $xE_*(G)y$ fails.

On the other hand, for $k\le k'$, if $2M_{r-1}\le k<2M_r$ and $2M_{r'-1}\le k'<2M_{r'}$, then
$$\begin{array}{ll}& d_H((S^C)^\omega(x),(S^C)^\omega(y))|_{[k,k']}\cr
&\le\sum_{p=2M_{r-1}}^{2M_{r'}-1}d_H(S^C(x(p)),S^C(y(p)))\cr
&=\sum_{q=r-1}^{r'-1}\sum_{i=0}^{m_q-1}d_H(S^C(g_{q,i}),S^C(h_qg_{q,i}))\cr
&<\sum_{q=r-1}^{r'-1}\sum_{i=0}^{m_q-1}2^{-(q+i+3)}<\sum_{q=r-1}^{r'-1}2^{-(q+2)}<2^{-r}.\end{array}$$
This gives $(S^C)^\omega(x)E_*(H)(S^C)^\omega(y)$. A contradiction!
\end{proof}

We conclude this subsection by the following theorem:

\begin{theorem}[Pre-rigid Theorem]\label{pre-rigid}
Let $G,H$ be two TSI Polish groups such that $H$ is locally compact. If $E(G)\le_BE(H)$, then there exist an open normal subgroups $G_c$ of $G$ and a continuous map $S:G_c\to H$ with $S(1_G)=1_H$ such that, for $x,y\in G_c^\omega$, if $\lim_nd_G(x(n),y(n))=0$, then
$$xE_*(G_c)y\iff S^\omega(x)E_*(H)S^\omega(y).$$

In particular, if $G=G_c$ and the interval $[0,1]$ embeds in $H$, then the converse is also true.
\end{theorem}

\begin{proof}
Suppose $E(G)\le_BE(H)$, then Lemma~\ref{S^C} gives the desired results.

Conversely, if $G=G_c$, from such an $S:G\to H$ we can define a continuous map $\theta:G^\omega\to H^\omega\times G^\omega$ as: for $x\in G^\omega$ and $n\in\omega$,
$$\theta(x(n))=(S(x(n)),x(n)).$$
It is clear that $S$ is a Borel reduction of $E_*(G)$ to $E_*(H)\times E(G;0)$. Since the interval $[0,1]$ embeds in $H$, following the arguments in the proof of Theorem~\ref{[0,1]}, we can see that $E(G)\le_BE(H)$.
\end{proof}

Under the assumption in Pre-rigid Theorem, the existence of Borel reduction is transformed into the existence of a well-behaved continuous mapping between topological groups.

\subsection{Applications on Lie groups}

A {\it Lie group} is a group which is also a smooth manifold such that the group operations are smooth functions. A Lie group is Polish iff it is separable iff it has only countablely many connected components. Let $G$ be a Lie group, we denote by $G_0$ the connected component of $1_G$, then $G_0$ is an open normal subgroup of $G$. For more details on Lie groups, we refer to~\cite{knapp,varadarajan}.

\begin{theorem}\label{Lie groups}
Let $G$ and $H$ be two separable TSI Lie groups. If $E(G)\le_BE(H)$, then there exists a continuous locally injective map $S:G_0\to H_0$.
\end{theorem}

\begin{proof}
If $E(G)\le_BE(H)$, then $E_*(G)\le_BE_*(H)$. The open normal subgroup $G_c$ of $G$ defined in the arguments before Lemma~\ref{S^C} is clopen, so $G_c$ contains the connected component $G_0$. By continuity of $S$, we have $S(G_0)\subseteq H_0$. Let $S^C$ be which appeared in Lemma~\ref{locally injective}. Note that any Lie group is locally compact and has no small subgroups (cf. \cite[Proposition 2.17]{HM07}). So the restriction of $S^C$ on $G_0$ is the desired continuous locally injective map.
\end{proof}

Let $f:\tilde X\to X$ be a continuous surjection between two topological spaces. We say $f$ is a {\it covering map} provided that, for any $x\in X$, there exists an open $V\ni x$ such that $f^{-1}(V)$ is the disjoint union of open subsets of $\tilde X$ each of which is mapped homeomorphically onto $V$ under $f$. Then $\tilde X$ is said to be a {\it covering space} of $X$.

\begin{theorem}\label{dimension}
Let $G$ and $H$ be two separable TSI Lie groups with $E(G)\le_BE(H)$. Then we have
\begin{enumerate}
\item[(1)] $\dim(G)\le\dim(H)$;
\item[(2)] if $\dim(G)=\dim(H)\ge 2$ and $G_0$ is compact, then $H_0$ is also compact and $G_0$ is a covering space of $H_0$.
\end{enumerate}
\end{theorem}

\begin{proof}
Clause (1) follows from Theorem~\ref{Lie groups} trivially.

If $\dim(G)=\dim(H)\ge 2$ and $G_0$ is compact, by Theorem~\ref{Lie groups}, there exists a continuous locally injective map $S:G_0\to H_0$. It follows from~\cite[Lemma 1]{MO} that $S$ is an open map, so $S(G_0)$ is an open subset of $H_0$. Since $G_0$ is compact, we have $S(G_0)$ is closed, and hence clopen in $H_0$. So $S(G_0)=H_0$, thus $H_0$ is compact. Moreover, since $S$ is an open continuous locally injective map, it is a locally homeomorphism. Then by~\cite[Propositon 1]{CP}, $S$ is a covering map.
\end{proof}

Clause (2) is also valid for $\dim(G)=\dim(H)=1$. We will show it later as a special case of that $G_0$ and $H_0$ are abelian.

\begin{remark}
Lemma 5.8 of~\cite{ding17} concerning finite dimensional Banach spaces is a special case of the preceding theorem. Indeed, it was this lemma from~\cite{ding17} that served as the inspiration for our definition of the equivalence relations induced by Polish groups.
\end{remark}

It is well known that all abelian Polish groups and compact Polish groups are TSI (cf.~\cite[Exercise 2.1.5]{gaobook}). On the other hand, a connected locally compact Polish group is TSI iff it is the product of a compact Polish group and a group $\R^n$ (cf.~\cite[Exercise 2.1.4]{gaobook} and~\cite[16.4.6]{dixmier}). Let $G$ be a connected TSI Lie group. Since any closed subgroup of a Lie group is also a Lie group (cf.~\cite[Theorem 1.4]{arvanitoyeorgos}), we can write $G=\R^n\times G_K$, where $G_K$ is a compact connected Lie group.

Now we focus on the case that $G_0,H_0$ are abelian. It is well known that every connected abelian Lie group is of the form $\mathbb{R}^c\times\mathbb T^e$ (cf.~\cite[Proposition 1.12]{arvanitoyeorgos}). We can write
$$G_0=\R^{c_0}\times\T^{e_0},\quad H_0=\R^{c_1}\times\T^{e_1}.$$
Suppose $E(G)\le_BE(H)$, we will show that $S^C:G_0\to H_0$ is a group homomorphism for each coset $C$ of $G_0$. Note that $\T=\{\exp(\i t):t\in\R\}$.

Given $g_0,g_1\in G_0$ and fix their expressions as
$$g_0=(\tau_0^0,\cdots,\tau_0^{c_0-1},\exp(\i\tau_0^{c_0}),\cdots,\exp(\i\tau_0^{c_0+e_0-1})),$$
$$g_1=(\tau_1^0,\cdots,\tau_1^{c_0-1},\exp(\i\tau_1^{c_0}),\cdots,\exp(\i\tau_1^{c_0+e_0-1})).$$
For $t\in[0,1]$, define $g(t)\in G_0$ such that
$${\rm proj}_i(g(t))=\left\{\begin{array}{ll}\tau_0^i+t(\tau_1^i-\tau_0^i), & 0\le i<c_0,\cr
\exp(\i(\tau_0^i+t(\tau_1^i-\tau_0^i))), & c_0\le i<c_0+e_0.\end{array}\right.$$
Then $g(0)=g_0$ and $g(1)=g_1$.

Fix an $i_0<c_1+e_1$. If $i_0<c_1$, we define $F_{i_0}:[0,1]\to\R$ as $F_{i_0}(t)={\rm proj}_{i_0}(S^C(g(t)))$. Note that $[0,1]$ is simply connected. If $c_1\le i_0<c_1+e_1$, by~\cite [Definition A2.6]{HM13}, we can also find a continuous function $F_{i_0}:[0,1]\to\R$ such that $\exp(\i F_{i_0}(t))={\rm proj}_{i_0}(S^C(g(t)))$.

\begin{lemma}\label{linear}
$F_{i_0}(t)=F_{i_0}(0)+t(F_{i_0}(1)-F_{i_0}(0))$ for $t\in[0,1]$.
\end{lemma}

\begin{proof}
Define $f:[0,1]\to\R$ as $f(t)=F_{i_0}(t)-F_{i_0}(0)-t(F_{i_0}(1)-F_{i_0}(0))$. Note that $f$ is continuous and $f(0)=f(1)=0$. We only need to prove that $f(t)=0$ for all $t\in(0,1)$.

If not, there exists $t_0\in(0,1)$ such that $f(t_0)\ne 0$. Without loss of generality, assume that $f(t_0)>0$. Let
$$\xi=\inf\{t\le t_0:f(t)\ge f(t_0)\},\quad \zeta=\sup\{t\ge t_0:f(t)\ge f(t_0)\}.$$
By the continuity, we have $f(\xi)=f(\zeta)=f(t_0)$. Again by the continuity, we can find $0<\xi_0<\xi_1<\xi_2<\cdots<\xi$ such that $f(\xi_l)=\frac{l+1}{l+2}f(t_0)$ for each $l\in\omega$. Let $\xi'=\lim_l\xi_l$, then $\xi'\le\xi$ and $f(\xi')=\lim_lf(\xi_l)=f(t_0)$. By the definition of $\xi$, we have $\xi=\xi'$. There exists $L\in\omega$ such that $\xi-\xi_L<1-\zeta$. Thus we can find $1>\zeta_0>\zeta_1>\zeta_2>\cdots>\zeta$ such that, for $l\ge L$, we have
$$\xi-\xi_l>\zeta_l-\zeta>\xi-\xi_{l+1}.$$
This gives $\zeta=\lim_l\zeta_l$.

For $p\in\omega$, we set
$$x(p)=\left\{\begin{array}{ll}g(\xi), & p=2l,\cr g(\zeta), & p=2l+1,\end{array}\right.\quad
y(p)=\left\{\begin{array}{ll}g(\xi_l), & p=2l,\cr g(\zeta_l), & p=2l+1.\end{array}\right.$$
It is clear that $\lim_pd_G(x(p),y(p))=0$, so Lemma~\ref{S^C} gives
$$xE_*(G)y\iff (S^C)^\omega(x)E_*(H)(S^C)^\omega(y).$$

From the alternating series test, the following series
$$(\xi-\xi_0)+(\zeta-\zeta_0)+\cdots+(\xi-\xi_l)+(\zeta-\zeta_l)+\cdots$$
is convergent. We also note that
$${\rm proj}_i(g(\xi)g(\xi_l)^{-1})=
\left\{\begin{array}{ll}(\xi-\xi_l)(\tau_1^i-\tau_0^i), & 0\le i<c_0,\cr
\exp(\i(\xi-\xi_l)(\tau_1^i-\tau_0^i)), & c_0\le i<c_0+e_0,\end{array}\right.$$
$${\rm proj}_i(g(\zeta)g(\zeta_l)^{-1})=
\left\{\begin{array}{ll}(\zeta-\zeta_l)(\tau_1^i-\tau_0^i), & 0\le i<c_0,\cr
\exp(\i(\zeta-\zeta_l)(\tau_1^i-\tau_0^i)), & c_0\le i<c_0+e_0.\end{array}\right.$$
It follows that, for all $i<c_0+e_0$, the following are convergent:
$$\sum_p{\rm proj}_i(x(p)y(p)^{-1})=\sum_l{\rm proj}_i(g(\xi)g(\zeta)g(\xi_l)^{-1}g(\zeta_l)^{-1})\quad(0\le i<c_0),$$
$$\prod_p{\rm proj}_i(x(p)y(p)^{-1})=\prod_l{\rm proj}_i(g(\xi)g(\zeta)g(\xi_l)^{-1}g(\zeta_l)^{-1})\quad(c_0\le i<c_0+e_0).$$
So we have $xE_*(G)y$.

On the other hand, we have
$$\sum_l(f(\xi)-f(\xi_l)+f(\zeta)-f(\zeta_l))\ge\sum_l(f(\xi)-f(\xi_l))=\sum_l\frac{f(t_0)}{l+2}=\infty.$$
Since
$$\begin{array}{ll}& F_{i_0}(\xi)-F_{i_0}(\xi_l)+F_{i_0}(\zeta)-F_{i_0}(\zeta_l)\cr
=&(f(\xi)-f(\xi_l)+f(\zeta)-f(\zeta_l))+(\xi-\xi_l+\zeta-\zeta_l)(F_{i_0}(1)-F_{i_0}(0)),\end{array}$$
and $\sum_l(\xi-\xi_l+\zeta-\zeta_l)$ converges, we have
$$\sum_l(F_{i_0}(\xi)-F_{i_0}(\xi_l)+F_{i_0}(\zeta)-F_{i_0}(\zeta_l))=\infty.$$
If $i_0<c_1$, we have
$$\sum_l{\rm proj}_{i_0}(S^C(g(\xi))S^C(g(\xi_l))^{-1}S^C(g(\zeta))S^C(g(\zeta_l))^{-1})=\infty;$$
if $c_1\le i_0<c_1+e_1$, since $(F_{i_0}(\xi)-F_{i_0}(\xi_l)+F_{i_0}(\zeta)-F_{i_0}(\zeta_l))\to 0$,
$$\begin{array}{ll} & \prod_l{\rm proj}_{i_0}(S^C(g(\xi))S^C(g(\xi_l))^{-1}S^C(g(\zeta))S^C(g(\zeta_l))^{-1})\cr
=& \prod_l\exp(\i(F_{i_0}(\xi)-F_{i_0}(\xi_l)+F_{i_0}(\zeta)-F_{i_0}(\zeta_l)))\mbox{ diverges}.\end{array}$$

Therefore, $\prod_pS^C(x(p))S^C(y(p))^{-1}$ diverges. In other words, we have $(S^C)^\omega(x)E_*(H)(S^C)^\omega(y)$ fails. A contradiction!
\end{proof}

\begin{lemma}\label{homo abelian}
$S^C:G_0\to H_0$ is a group homomorphism.
\end{lemma}

\begin{proof}
For any $i_0<c_1+e_1$, Lemma~\ref{linear} gives $F_{i_0}(1/2)=F_{i_0}(0)+(F_{i_0}(1)-F_{i_0}(0))/2$, i.e. $F_{i_0}(0)+F_{i_0}(1)=2F_{i_0}(1/2)$. It follows that
$$S^C(g(0))S^C(g(1))=S^C(g(1/2))^2.$$

For $t\in[0,1]$, define $g^*(t)\in G_0$ as
$${\rm proj}_i(g^*(t))=\left\{\begin{array}{ll}t(\tau_0^i+\tau_1^i), & 0\le i<c_0,\cr
\exp(\i(t(\tau_0^i+\tau_1^i))), & c_0\le i<c_0+e_0.\end{array}\right.$$
Similarly, we have
$$S^C(g^*(0))S^C(g^*(1))=S^C(g^*(1/2))^2.$$

From the definition of $S^C$, we have $S^C(1_G)=1_H$. Note that $g(0)=g_0$ and $g(1)=g_1$, also note that $g^*(0)=1_G,g^*(1)=g_0g_1$, and $g^*(1/2)=g(1/2)$. So we have
$$S^C(g_0g_1)=S^C(g^*(1))=S^C(g^*(1/2))^2=S^C(g(1/2))^2=S^C(g_0)S^C(g_1).$$
Therefore, $S^C$ is a group homomorphism.
\end{proof}

\begin{theorem}
$E(\R^{c_0}\times\T^{e_0})\le_BE(\R^{c_1}\times\T^{e_1})$ iff $e_0\le e_1$ and $c_0+e_0\le c_1+e_1$.
\end{theorem}

\begin{proof}
The $(\Leftarrow)$ part follows form Theorem~\ref{[0,1]} and Proposition~\ref{subgroup}.

Now suppose $E(\R^{c_0}\times\T^{e_0})\le_BE(\R^{c_1}\times\T^{e_1})$. Theorem~\ref{dimension} gives $c_0+e_0\le c_1+e_1$. Since $E(\T^{e_0})\le_BE(\R^{c_0}\times\T^{e_0})\le_BE(\R^{c_1}\times\T^{e_1})$, by Theorem~\ref{Lie groups} and Lemma~\ref{homo abelian}, there exists a continuous locally injective homomorphism $S:\T^{e_0}\to\R^{c_1}\times\T^{e_1}$. Then there are continuous homomorphisms $R:\T^{e_0}\to\R^{c_1}$ and $T:\T^{e_0}\to\T^{e_1}$ such that $S(g)=(R(g),T(g))$ for $g\in\T^{e_0}$. It is trivial to see that $R(g)=0$ for all $g\in\T^{e_0}$. So $T$ should be locally injective, and hence $e_0\le e_1$.
\end{proof}

Generalizing the results above, the authors~\cite{DZtsi} proved a Rigid Theorem concerning TSI Lie groups as following:

\begin{thm}[{\cite[Theorem 1.8]{DZtsi}}]
Let $G,H$ be two separable TSI Lie groups such that $G$ is connected. Then $E(G)\leq_B E(H)$ iff there exists a continuous locally injective homomorphism $S:G\rightarrow H$.
\end{thm}

Next we provide an example of a disconnected Lie group. Let $\Z_2=\{0,1\}$ and denote by $\Z_2\ltimes\R$ the semi product $\Z_2\ltimes_\phi\R$ with $\phi(0)={\rm id_\R}$ and $\phi(1)(t)=-t$ for $t\in\R$. Then we have $(1,0)(0,t)(1,0)=(0,-t)$ for $t\in\R$. In the following, we write $v=(1,0)$ for brevity. Note that $v^{-1}=v$ and $v(0,t)v^{-1}=(0,-t)$.

\begin{lemma}\label{rtimes}
$E(\Z_2\ltimes\R)\not\le_BE(\T)$.
\end{lemma}

\begin{proof}
We use notions defined in the arguments before Lemma~\ref{S^C}. We can assume $v\in F_0$. In order to write succinctly, we identify $\{0\}\times\R$ with $\R$.

Assume for contradiction that $E(\Z_2\ltimes\R)\le_BE(\T)$. Let $G=\Z_2\ltimes\R$ and $H=\T$. Since $\T$ is compact, we can put $H_c=W_c=\T$. Then we can put $G_c=V_c$ to be either $\R$ or $G$.

First, for $G_c=\R$, there are two cosets $\R$ and $\R v$. By Lemma~\ref{homo abelian}, $S^\R,S^{\R v}:\R\to\T$ are continuous locally injective homomorphisms. Thus there exist $a,b\in\R\setminus\{0\}$ such that, for $t\in\R$,
$$S^\R((0,t))=\exp(\i at),\quad S^{\R v}((0,t))=\exp(\i bt).$$

Secondly, for $G_c=G$, there is only one coset $G$. Since $H=\T$ is abelian and $v^2=1_G$, we have
$$S^{\R v}_{n_{2j}}(g)=S_{n_{2j}}(gv)S_{n_{2j}}(v)^{-1},\quad S^G_{n_{2j}}(h)=S_{n_{2j}}(h)S_{n_{2j}}(1_G)^{-1}$$
for $g=(0,t)\in\R\cap F_{n_{2j}-2},~v\in F_{n_{2j}-2}$ and any $h\in F_{n_{2j}-1}$, so
$$S^{\R v}_{n_{2j}}((0,t))=S_{n_{2j}}((0,t)v)S_{n_{2j}}(v)^{-1}=S^G_{n_{2j}}((0,t)v)S^G_{n_{2j}}(v)^{-1}.$$
Assume that the ultrafilter $\frak A\ni\{2j:j\in\omega\}$, then we have
$$S^{\R v}((0,t))=S^G((0,t)v)S^G(v)^{-1}.$$
So $S^G((0,t)v)=\exp(\i bt)S^G(v)$ for $t\in\R$.

Now define $x,y\in G^\omega$ as $x(p)=v$ and $y(p)=(0,\frac{1}{p+1})v$ for $p\in\omega$. Then we have $\lim_pd_G(x(p),y(p))=0$, and hence by Lemma~\ref{S^C},
$$xE_*(G)y\iff (S^G)^\omega(x)E_*(H)(S^G)^\omega(y).$$

For $p\in\omega$, we have $x(0)\cdots x(p)=v^{p+1}$ and
$$\begin{array}{ll}y(0)\cdots y(p)&=(0,1)v\cdots(0,\frac{1}{p+1})v\cr
&=(0,1)(v(0,\frac{1}{2})v^{-1})\cdots(v^p(0,\frac{1}{p+1})v^{-p})v^{p+1}\cr
&=(0,1+\frac{-1}{2}\cdots+\frac{(-1)^p}{p+1})v^{p+1}.\end{array}$$
Since $\sum_p\frac{(-1)^p}{p+1}$ converges, we have $xE_*(G)y$.

On the other hand, we have $S^G(x(0))\cdots S^G(x(p))=S^G(v)^{p+1}$ and
$$\begin{array}{ll}S^G(y(0))\cdots S^G(y(p))&=\exp(\i b)S^G(v)\cdots\exp(\frac{\i b}{p+1})S^G(v)\cr
&=\exp(\i b)\exp(\frac{\i b}{2})\cdots\exp(\frac{\i b}{p+1})S^G(v)^{p+1}.\end{array}$$
Since $\sum_p\frac{1}{p+1}=\infty$ and $\frac{1}{p+1}\to 0$, we have $\prod_p\exp(\frac{\i b}{p+1})$ diverges, and hence $(S^G)^\omega(x)E_*(H)(S^G)^\omega(y)$ fails. A contradiction!
\end{proof}

We also denote by $\Z_2\ltimes\T$ the semi product $\Z_2\ltimes_\phi\T$ with $\phi(0)={\rm id_\T}$ and $\phi(1)(\exp(\i t))=\exp(-\i t)$ for $t\in[0,2\pi)$.

\begin{corollary}\label{RT}
\begin{enumerate}
\item[(1)] $E(\T)$ and $E(\Z_2\ltimes\R)$ are Borel incomparable;
\item[(2)] $E(\R)<_BE(\T)<_BE(\Z_2\ltimes\T)$;
\item[(3)] $E(\R)<_BE(\Z_2\ltimes\R)<_BE(\Z_2\ltimes\T)$.
\end{enumerate}
\end{corollary}

\begin{proof}
Clause (1) follows from lemmas~\ref{homo abelian} and~\ref{rtimes}. By Proposition~\ref{subgroup} and Theorem~\ref{[0,1]}, we have
$E(\R)\le_BE(\T)\le_BE(\Z_2\ltimes\T)$ and $E(\R)\le_BE(\Z_2\ltimes\R)\le_BE(\Z_2\ltimes\T)$. Then (2) and (3) follow from (1).
\end{proof}

\begin{corollary}\label{innR}
Let $G$ be a separable TSI Lie group with $G_0=\R$. Define $\iota_\R:\R\to\R$ as $\iota_\R(t)=-t$. Then we have
\begin{enumerate}
\item[(1)] $E(G)\sim_BE(\R)\iff{\rm Inn}_G(G_0)=\{{\rm id}_\R\}$;
\item[(2)] $E(G)\sim_BE(\Z_2\ltimes\R)\iff{\rm Inn}_G(G_0)=\{{\rm id}_\R,\iota_\R\}$.
\end{enumerate}
\end{corollary}

\begin{proof}
Suppose that there are $u\in G$ and $b\in\R$ with $b\ne\pm 1$ such that $utu^{-1}=bt$ for all $t\in\R=G_0$. Without loss of generality, we can assume that $|b|<1$. Then $u^ntu^{-n}=b^nt\to 0$, contradicting that $G$ is TSI. So ${\rm Inn}_G(G_0)\subseteq\{{\rm id}_\R,\iota_\R\}$.

By Theorem~\ref{inn} and $E(\Z_2\ltimes\R)\not\le_BE(\R)$, we complete the proof.
\end{proof}

By similar arguments, we can prove the following corollary.

\begin{corollary}\label{innT}
Let $G$ be a separable TSI Lie group with $G_0=\T$. Define $\iota_\T:\T\to\T$ as $\iota_\T(\exp(\i t))=\exp(-\i t)$. Then we have
\begin{enumerate}
\item[(1)] $E(G)\sim_BE(\T)\iff{\rm Inn}_G(G_0)=\{{\rm id}_\T\}$;
\item[(2)] $E(G)\sim_BE(\Z_2\ltimes\T)\iff{\rm Inn}_G(G_0)=\{{\rm id}_\T,\iota_\T\}$.
\end{enumerate}
\end{corollary}

\subsection{Applications on $p$-adic solenoids}

Let $p\ge 2$ be a natural number. Recall that the {\it $p$-adic solenoid} $\T_p$ is the closed subgroup of $\T^\omega$ as follows:
$$\T_p=\{(g_l)\in\T^\omega:\forall l\,(g_l=g_{l+1}^p)\}.$$
Then $\T_p$ is a compact connected abelian Polish group, but not arcwise connected. The arc component of $1_{\T_p}$ in $\T_p$ is
$$R_p=\{(\exp(\i t/p^l)):t\in\R\}$$
(cf.~\cite[Exercise E1.11.(iv)]{HM13}). For more on $p$-adic solenoids, we refer to~\cite{HM13}.

Given $p\ge 2$, let $N(p)$ be the set of all prime divisors of $p$.

\begin{theorem}
Let $p,q\ge 2$ be natural numbers, then we have
\begin{enumerate}
\item[(1)] $E(\T_p)\le_BE(\T_q)\iff N(p)\supseteq N(q)$;
\item[(2)] $E(\R)<_BE(\T_p)<_BE(\T)$.
\end{enumerate}
\end{theorem}

\begin{proof}
For $p\ge 2$, let $\phi_p:\T_p\to\T$ as $\phi_p((g_l))=g_0$ for $(g_l)\in\T_p$. Then $\phi_p$ is a closed continuous surjective group homomorphism, so $\T\cong\T_p/\ker(\phi_p)$. Let $\psi_p:\R\to\T_p$ as $\psi_p(t)=(\exp(\i t/p^l))$ for $t\in\R$. Then $\psi_p$ is a continuous injective group homomorphism with $R_p=\psi_p(\R)$.

(1) First, if $q|p$, there is $d\in\N$ with $p=qd$. Let $\phi_{p,q}:\T_p\to\T_q$ as $\phi_{p,q}((g_l))=(g_l^{d^l})$ for $(g_l)\in\T_p$. Then $\phi_{p,q}$ is a closed continuous surjective group homomorphism with
$$\ker(\phi_{p,q})=\T_p\cap\{(g_l)\in\T^\omega:\forall l\,\exists k<d^l\,(g_l=\exp(2k\pi\i/d^l))\}.$$
It is clear that $\ker(\phi_{p,q})$ is non-archimedean and $\T_q\cong\T_p/\ker(\phi_{p,q})$. By Theorem~\ref{[0,1]}, we have $E(\T_p)\le_BE(\T_q)$.

Secondly, if $N(p)\supseteq N(q)$, there is $e\in\N$ with $q|p^e$. Let $\phi:\T_p\to\T_{p^e}$ as $\phi((g_l))=(g_{el})$ for $(g_l)\in\T_p$. Then $\phi$ is a topological group isomorphism. So $\T_p\cong\T_{p^e}$, and hence $E(\T_p)\sim_BE(\T_{p^e})\le_BE(\T_q)$.

On the other hand, if $E(\T_p)\le_BE(\T_q)$, since $\T_q$ is compact, by Lemma~\ref{S^C}, there exists a continuous map $S:\T_p\to\T_q$ with $S(1_{\T_p})=1_{\T_q}$ such that, for $x,y\in\T_p^\omega$, if $\lim_id_{\T_p}(x(i),y(i))=0$, then
$$xE_*(\T_p)y\iff S^\omega(x)E_*(\T_q)S^\omega(y).$$
Moreover, since $\phi_q$ is a continuous homomorphism, we have
$$S^\omega(x)E_*(\T_q)S^\omega(y)\Rightarrow\phi_q^\omega(S^\omega(x))E_*(\T)\phi_q^\omega(S^\omega(y)).$$

Given $\tau\in\R$, without loss of generality, assume that $\tau>0$. For $0\le t\le\tau$, $\psi_p(t)$ is in the arc component $R_p$ of $1_{\T_p}$, thus $S(\psi_p(t))$ is in the arc component $R_q$ of $1_{\T_q}$. Let $F:[0,\tau]\to\R$ be a continuous function such that $\exp(\i F(t))=\phi_q(S(\psi_p(t))))$ with $F(0)=0$. Following the arguments in the proof of Lemma~\ref{linear}, we have $F(t)=tF(\tau)/\tau$. Due to the arbitrariness of $\tau$, there is $a\in\R$ such that $a=F(t)/t$ for all $t\neq 0$, so $S(\psi_p(t))=\psi_q(at)$. It is clear that $a\ne 0$.

Note that $\lim_m\psi_p(p^m)=1_{\T_p}$. By the continuity of $S$, $\lim_m\psi_q(ap^m)=1_{\T_q}$. For each $l\in\omega$, we have $\lim_m\exp(\i ap^m/q^l)=1$ and hence $\frac{ap^m}{2\pi q^l}\in\Z$ for $m$ large enough. Thus $N(p)\supseteq N(q)$.

(2) By (1) we only need to show $E(\R)\le_BE(\T_p)\le_BE(\T)$. Note that
$$\ker(\phi_p)=\T_p\cap\{(g_l)\in\T^\omega:\forall l\,\exists k<p^l\,(g_l=\exp(2k\pi\i/p^l))\}.$$
It is clear that $\ker(\phi_p)$ is non-archimedean. Thus by Theorem~\ref{[0,1]}, we have $E(\T_p)\le_BE(\T)$.

Note that $\phi_p(\psi_p(\R))=\T$ and $\ker(\phi_p\circ\psi_p)=\{2k\pi:k\in\Z\}$ is discrete. Since $[0,1]$ embeds into $\T_p$, by Theorem~\ref{[0,1]}, we have $E(\R)\le_BE(\T_p)$.
\end{proof}

\begin{remark}
The authors also generalized the above theorem to $P$-adic solenoids, where $P$ is a sequence of natural numbers $\ge 2$ (see~\cite[Theorem 3.2]{DZ}). Furthermore, the Borel reducibility among $E(G)$'s between $E(\R)$ and $E(\T)$ are extremely complicated that the partial ordered set $P(\omega)/{\rm Fin}$ embeds into them (see~\cite[Theorem 3.6]{DZ}).
\end{remark}

\section{Further remarks}

Although we have already presented many results concerning equivalence relations $E(G)$'s, more interesting questions remaind unanswered.
For instance, Theorem~\ref{trivial}(5) gives $E(S_\infty)\sim_B =^+$. And we already obtained a sequence of CLI non-archimedean Polish groups $(G_n)$ so that
$$E_0\sim_BE(G_0)<_B\cdots<_BE(G_n)<_BE(G_{n+1})<_B\cdots<_B E(S_\infty)\sim_B=^+.$$
However, we do not know:

\begin{question}
Is there a sequence of non-archimedean Polish groups $(G_n)_{n\in\omega}$ such that these $E(G_n)$'s are pairwise Borel incomparable?
\end{question}

In section 5, we use the notion of l.m.-unbalanced groups to prove some $E(G)$ are not Borel reducible to $E(H)$ for any TSI Polish group $H$. As mentioned in the first paragraph of subsection 5.3, all nilpotent topological groups are distal. And by Theorem~\ref{lc-nondistal}, we know all nilpotent locally compact Polish groups are not l.m.-unbalanced. For instance, let $T^u(3)$ be the group of all real matrices as $\left(\begin{array}{lcr}1 & a & b\cr 0 & 1 & c\cr 0 & 0 & 1\end{array}\right)$. It is clear that $T^u(3)$ is nilpotent locally compact, but not TSI.
It follows from~\cite[Theorem 1.2]{DZtsi} that $E(T^u(3))\nleq_B E(H)$ for any TSI Polish group $H$. In particular, we have $E(T^u(3))\nleq_B E(\mathbb R^3)$.
So far, we do not know:

\begin{question}
Does $E(\mathbb R^3)\leq_B E(T^u(3))$?
\end{question}

As mentioned at the end of section 5, we do not know any example of $\omega$-l.m.-unbalanced CLI Polish group. So perhaps the most noteworthy problem is:

\begin{question}
Does $\alpha$-l.m.-unbalanced CLI Polish groups exist for each $\alpha<\omega_1$? If not, can we find more sophisticated tools to form a hierarchy of length $\omega_1$ for all CLI Polish groups under Borel reducibility among equivalence relations induced by them?
\end{question}

Some questions asked in earlier versions of this article have been answered in~\cite{DZ,DZtsi}. For instance, Theorem 1.8 of~\cite{DZtsi} is a positive answer for the following question.

\begin{question}\label{Qconnected}
Let $G,H$ be two connected TSI Lie groups, does $E(G)\le_BE(H)$ iff there exists a continuous locally injective homomorphism from $G$ to $H$?
\end{question}

\begin{question}\label{Qdisconnected}
Let $G,H$ be two separable TSI Lie groups such that their identity component $G_0$ and $H_0$ are topologically isomorphic, does $E(G)\le_BE(H)$ iff there exists an topological isomorphism $\phi:G_0\to H_0$ such that $\phi{\rm Inn}_G(G_0)\phi^{-1}\subseteq{\rm Inn}_H(H_0)$?
\end{question}

Theorem~\ref{homo abelian} confirms the special case of $G,H$ are abelian for Question~\ref{Qconnected}. And corollaries~\ref{innR} and~\ref{innT} confirm the very special case of $G_0$ and $H_0$ are topologically isomorphic to either $\R$ or $\T$ for Question~\ref{Qdisconnected}.

The authors~\cite{DZtsi} considered equivalence relations induced by separable strongly NSS Fr{\'e}chet spaces.
So far, we know almost nothing about Borel reducibility among equivalence relations induced by non-strongly NSS  infinitely dimensional Polish groups.

\subsection*{Acknowledgements}
We thank Su Gao, Feng Li, Yiming Sun, Ruijun Wang and Xu Wang for helpful discussions in seminars. We are grateful to Xiangjun Wang, Shizhuo Yu, Hui Zhang and Fuhai Zhu for their advices on algebraic topology and Lie groups. Special thanks are due to Cheng Peng for carefully reading and suggestions on proofs of some results.

\end{document}